\documentclass[reqno]{amsart}
        \newif\ifpdf 
        \ifx\pdfoutput\undefined 
        \pdffalse % we are not running pdflatex 
        \else 
        \pdfoutput=1 % we are running pdflatex 
        \pdftrue 
        \fi

\usepackage{latexsym}
\usepackage{amstext}
\usepackage {amsmath}
\usepackage {amsfonts}
\usepackage {amssymb}
\usepackage {amsthm}
\usepackage {bbm}
\usepackage{enumitem}%\usepackage{enumerate}
\usepackage{array}
\usepackage{comment}
\usepackage{xspace}
\DeclareMathOperator{\dist}{dist}
\DeclareMathOperator{\graph}{graph}

\DeclareMathOperator{\spt}{supp}

\def\loc{{\mathrm{loc}}}
\DeclareMathOperator{\Lp}{L}
\DeclareMathOperator{\BV}{BV}

\DeclareMathOperator{\trace}{tr}

\newcommand{\eps}{\varepsilon}
\newcommand{\Chi}{\mathcal{X}}
\newcommand{\R}{\ensuremath{\mathbb{R}}}
\newcommand{\Rn}{\ensuremath{\mathbb{R}^n}}
\newcommand{\N}{\ensuremath{\mathbb{N}}\ }
\newcommand{\LL}{\ensuremath{\mathcal{L}}}

\newcommand{\Ha}{\ensuremath{\mathcal{H}}}

\newcommand{\ddist}{\ensuremath{\mathsf{d}}}

\def\step#1{\\[2pt]\textsc{Step#1.}}

\newcommand{\fplus}{F_+}
\newcommand{\fmin}{F_-}
\newcommand{\Sym}{\mathcal{S}}
\newcommand{\smallc}{c_1}
\newcommand{\smalleta}{\omega}
\newcommand{\smallbdry}{\omega}
\newcommand{\MK}{\mathsf{H}}
\newcounter{counter-liste}
\newcounter{counter-liste2}

  {\ \begin{list}{{(\arabic{counter-liste})}\hfill}%
  {\topsep2mm\itemindent1ex\leftmargin0cm\usecounter{counter-liste}}
  }%
  {\end{list}}
\newenvironment{liste(a)}%
  {\ \begin{list}{{(\alph{counter-liste})}\hfill}%
  {\topsep2mm\itemindent1ex\leftmargin0cm\usecounter{counter-liste}}
  }%
  {\end{list}}
\def\R{\mathbb R}

\def\BV{\mathrm{BV}}

\theoremstyle{plain}
\numberwithin{equation}{section}
\newtheorem{lemma}{Lemma}[section]
\newtheorem{theorem}[lemma]{Theorem}
\newtheorem{proposition}[lemma]{Proposition}
\newtheorem{definition}[lemma]{Definition}
\newtheorem{assumption}[lemma]{Assumption}

\theoremstyle{definition}
\newtheorem{remark}[lemma]{Remark}
\begin{document}
\title[Convergence to the Gibbs--Thomson law]{Convergence of
  phase--field approximations to the Gibbs--Thomson law}  
\author{Matthias R{\"o}ger}
\address{Matthias R\"{o}ger, Max Planck Institute for Mathematics in the
  Sciences, Inselstr. 22, D-04103 Leipzig}

\author{Yoshihiro Tonegawa}
\address{Yoshihiro Tonegawa, Department of Mathematics,
Hokkaido University Sapporo, 060-0810 Japan}

\email{roeger@mis.mpg.de, tonegawa@math.sci.hokudai.ac.jp}



\date{\today}

\begin{abstract}
We prove the convergence of phase-field approximations of the
Gibbs--Thomson law. This establishes a relation between the first variation
of the Van-der-Waals--Cahn--Hilliard energy and the first variation
of the area functional. We allow for folding of diffuse interfaces in
the limit and the occurrence of higher-multiplicities of the limit energy
measures. We show that the multiplicity does not affect the
Gibbs--Thomson law and that the mean curvature vanishes where diffuse
interfaces have collided.

We apply our results to prove the convergence of
stationary points of the Cahn--Hilliard equation to constant mean
curvature surfaces and the convergence of
stationary points of an energy functional that was proposed by
Ohta--Kawasaki as a model for micro-phase separation in
block-copolymers. 
\end{abstract}
\maketitle
%\setcounter{tocdepth}{1}
%\tableofcontents
%==============================================
% intro
%==============================================
\section{Introduction}\label{sec:intro}
Phase separation is a common phenomenon in many
areas of the sciences.
Alloys studied in material sciences, melting and solidification
processes, or block-copolymers investigated in physical chemistry, they
all show the coexistence of two or more phases, separated by thin
\emph{transition layers}. The main approaches to describe phase
transitions are on the one hand \emph{sharp
  interface models} and  on the other hand \emph{diffuse interface
  models}, also referred to as `phase field'  or `Ginzburg-Landau'
models. The relation between both kinds of models remains an
outstanding question.  
Rigorous passages to the sharp interface limit are often
difficult and generalized formulations for the limit problems are
necessary to obtain the convergence of diffuse approximations.
However, care has to be taken that solutions satisfy the equations in
a reasonably strong sense. 

The goal of the present paper is to prove the convergence of diffuse
approximations of the so-called \emph{Gibbs--Thomson law}, which states
that the mean curvature of the phase boundary is given as the trace of a
function in the bulk. Our result relates the first
variation of the Van-der-Waal-Cahn--Hilliard energy, which is the common
root of most phase field models, to the first variation of the area
functional. To the best of our knowledge, we give
the first satisfactory solution in the case that diffuse interfaces
collapse or cancel each other in the limit. 

Before stating the main result we describe the setting and background of the
problem.
\subsection{Phase fields, sharp interfaces, and the Gibbs--Thomson
  law} 
The diffuse interface approach is based on a free energy that acts on
\emph{smooth phase fields} and that was proposed by Van-der-Waals
\cite{VdW} and later Cahn--Hilliard \cite{CH}.
In a normalized form this energy is given by
\begin{gather}
  E_{\varepsilon}(u)\,:=\,\int_{\Omega}\Big(\frac{\varepsilon}{2}|\nabla u|^2 
  +\frac{1}{\varepsilon}W(u)\Big) dx, \label{eq:def-E-intro}
\end{gather}
where $\eps>0$ is a small parameter and $W$ is a nonnegative `double-well
potential' with value zero if and only if $u=\pm 1$. 
Domains where $u\approx 1$ or $u\approx -1$ represent two coexisting
phases, separated by \emph{diffuse interfaces}.
Formal arguments show that $E_\eps$ favors transition layers with a
thickness of order $\varepsilon$. Hence, as $\eps$ tends to zero the diffuse
interfaces become sharp.  

One naturally associated quantity to the Cahn--Hilliard energy is its
$L^2$-functional derivative, which often corresponds to
the \emph{chemical potential}, 
\begin{gather}
   f_\eps \,=\, -\varepsilon\Delta u+\frac{1}{\varepsilon}W'(u). \label{eq:chem-pot} 
\end{gather}
In many applications $f_\eps$ is given by means of other quantities and
a certain control on $f_\eps$ is available.
The corresponding functional derivative of the surface area functional,
evaluated at a smooth compact hypersurface $\Sigma$, is given by the
mean curvature of $\Sigma$ and \eqref{eq:chem-pot} formally
corresponds to the following equation, in solidification processes known
as \emph{Gibbs--Thomson law} (and we will adopt this term throughout the
paper), 
\begin{gather}
  H \,=\,\sigma f, \label{eq:GT-H}
\end{gather}
with a surface tension coefficient $\sigma>0$. The Gibbs--Thomson law
relates the local geometry of the phase boundary to a function
$f:\Omega\to\R$ in the bulk, for example the temperature or the chemical
potential.  
\subsection{Main results.}
Let us first state our main result in a
concise form (we will prove a slightly stronger statement, given in Theorem
\ref{the:claim}). 
\begin{theorem}
\label{the:main-intro}
Suppose $p>n$ and let sequences of functions
$\{u_{\varepsilon}\}_{\eps>0}\subset W^{3,p}(\Omega)$ and functions
$\{f_\eps\}_{\eps>0}\subset W^{1,p}(\Omega)$ be given such that
\eqref{eq:chem-pot} holds and such that
\begin{alignat}{2}
  E_{\varepsilon}(u_{\varepsilon})\,&\leq\, \Lambda, \label{ass:bound-E-intro}\\
  u_\eps\,&\to\,u \quad&&\text{ in }L^1(\Omega), \label{ass:u-intro}\\
  f_\eps\,&\to\,f \quad&&\text{ weakly in }W^{1,p}(\Omega) \label{ass:f-intro}
\end{alignat}
as $\eps\to 0$.
Then $u$ is of bounded variations and takes only values in $\{-1,1\}$. Moreover
there exists a unique generalized mean curvature $H$ of the `phase
boundary' $\Sigma:=\Omega\cap\partial^{*}\{u=1\}$ and
\begin{gather}
  \sigma H \,=\, f\quad\text{ holds }\Ha^{n-1}-\text{almost everywhere on
  }\Sigma. \label{eq:intro-GT}
\end{gather}
Here $\sigma=\int_{-1}^1\sqrt{W(s)/2}\, ds$ is the surface tension
coefficient and the sign of $H$ is chosen positive for spherical
$\{u=1\}$.
\end{theorem}
This Theorem uses generalized formulations for the phase boundary and
the mean curvature. The notion of mean curvature is based on a
measure-theoretic approach and was introduced by the
first author in \cite{Ro04}. We refer to the appendix \ref{app:defs}
for the exact definition.
\begin{remark}\label{rem:main-intro}
For a sequence $(u_\eps)_{\eps>0}$ that satisfies the uniform energy bound
\eqref{ass:bound-E-intro} and a sequence $(f_\eps)_{\eps>0}$ that is
uniformly bounded in $W^{1,p}(\Omega)$ there exists a subsequence
$\eps\to 0$ such that \eqref{ass:u-intro}, \eqref{ass:f-intro} hold.
Besides these uniform bounds no other conditions, such as energy
minimality, are required to apply Theorem \ref{the:main-intro}. 
For this reason the result is
relevant to 
a large class of stationary and time-dependent problems. In section
\ref{sec:applications} we use our results to characterize the limit of
stationary points of the Cahn-Hilliard functional \eqref{eq:def-E-intro} and
to prove the convergence of stationary points in a model for
block-copolymers.

The assumption \eqref{ass:f-intro} on the chemical potentials $f_\eps$ is
still restrictive. We conjecture that the (weak) convergence of $f_\eps$
in $W^{1,p}(\Omega)$ with $p>n/2$ would suffice to conclude
\eqref{eq:intro-GT}.
However, our techniques yet require the continuity of $f$, which is
ensured only if $p>n$. For the Cahn--Hilliard equation for
example the natural regularity of the chemical potential is
$W^{1,2}(\Omega)$ in space. Hence, our result does not apply in this case.
\end{remark} 
\subsection{Related results and main techniques}
\label{subsec:related}
Since the fundamental work of Modica and Mortola \cite{MM,Mod} on the
convergence of $E_{\varepsilon}$ to the area functional the relation
between their first variations has drawn attention.
Modica \cite{Mod} and Sternberg \cite{Ste}
proved that minimizers of $E_\eps$ under a volume constraint
converge to area-minimizing hypersurfaces with an integral
constraint. Luckhaus--Modica \cite{LM} then showed that the
Lagrange-multipliers associated with the volume constraint converge to
the constant mean curvature of the limiting hypersurface. Ilmanen
\cite{Ilm} considered the corresponding $L^2$-gradient flows and proved
the convergence of the Allen--Cahn equation to the mean-curvature flow,
in the varifold formulation of Brakke \cite{Bra}. Convergence of various
other phase field problems to the corresponding sharp interface models
have been shown either formally or rigorously
\cite{Che92,ABC,Peg,CC98,Lu90,Che,So95}, sometimes in
quite involved weak formulations.

The second author considered, partly in joint-work with Hutchinson, the
convergence of diffuse interface approximations of the Gibbs--Thomson
law, under different assumptions on the chemical potential
\cite{HT,To02,To05}.
However, the Gibbs--Thomson relation is only verified in an
(in some respect unsatisfactory) multiplicity-dependent formulation,
see \eqref{eq:GT-mult} and the discussion below.

Sch\"atzle \cite{Sch01} considered a sequence of hypersurfaces with mean
curvature given by a Sobolev function in the ambient space and obtained
that the Gibbs--Thomson law holds in the limit in a rather clean
varifold formulation. In \cite{Sch01} the chemical potentials $f_\eps$
in \eqref{ass:f-intro} need only  to converge in a Sobolev space
$f\in W^{1,p}(\Omega)$ with $p>n/2$, \emph{c.f.} Remark
\ref{rem:main-intro}.

Geometric Measure Theory provides suitable 
generalized formulations in spaces that allow for the compactness of
approximations. 
Luckhaus--Modica \cite{LM} and Luckhaus--Sturzenhecker
\cite{LS} introduced a weak formulation of the Gibbs--Thomson law
\eqref{eq:GT-H} for characteristic functions of bounded
variation. This formulation is rather natural and has the advantage of
being based directly on the phase function. However,
justifying the Gibbs--Thomson law in the limit of approximations
requires the additional assumption that no cancellation of (diffuse)
interfaces occurs. Unfortunately, this property does in general not hold
\cite{Sch97}.

To master such cancellations Ilmanen \cite{Ilm} used a \emph{varifold-approach}.
He considered the limit of the \emph{diffuse surface-area measures}
(\emph{energy measures})
\begin{gather}
  \mu_\eps \,:=\,  \Big(\frac{\varepsilon}{2}|\nabla u_\eps|^2 
  +\frac{1}{\varepsilon}W(u_\eps)\Big)\,\LL^n.
\end{gather}
The idea behind is that this limit makes information visible
that is lost in the limit of the phase fields: Where cancellation of the
approximate phase boundaries occurs the limit $\mu$ of the measures
$\mu_\eps$ carries a higher multiplicity. The support of $\mu$
eventually extends the limit phase boundary by \emph{hidden
  boundaries}. Showing that the limit 
measure is in fact given as a integer-rectifiable varifold with a weak
mean curvature vector, the Gibbs--Thomson law can be verified in
a varifold formulation.

This strategy was used for various problems by Chen \cite{Che}, Soner
\cite{Son}, Hutchinson--Tonegawa \cite{HT}, Tonegawa
\cite{To02,To05}, and others. However, in none of these papers the
problem of higher multiplicity was completely solved.
Typically the convergence of the diffuse phase fields and the diffuse
surface-area measures $\mu_\eps$ is shown and the rectifiability of
the limit $\mu$ as well as the  existence
of a weak mean curvature $H_\mu$ is obtained. Still, the
Gibbs--Thomson law holds only in a multiplicity-dependent formulation
\begin{gather}
  H_\mu\,=\, (N\sigma)^{-1}f \quad \Ha^{n-1}\text{-almost everywhere on }
  \partial^*\{u=1\}, \label{eq:GT-mult}
\end{gather}
where $N$ is the density function of the measure $\mu$.
This formulation is for two reason unsatisfactory: First the
Gibbs--Thomson law should be satisfied by the \emph{phase boundary}
rather than by the (in view of the applications) `obscure' measure $\mu$.
Secondly, the density function $N$ should not affect the
Gibbs--Thomson law. To prove the full results \eqref{eq:GT-mult} has to
be complemented by
\begin{gather}\label{eq:f-H-0}
  f\,=\,H_\mu\,=\,0\quad \Ha^{n-1}\text{-almost everywhere on }
  \partial^*\{u=1\}\cap\{N>1\}.
\end{gather}
In recent years progress has been made on this issue. Sch\"atzle
\cite{Sch01} proves  the
Gibbs--Thomson law in the limit of an approximation by
hypersurfaces: There the weak mean curvature $H_\mu$ of the limit
measure $\mu$ satisfies \eqref{eq:GT-H} and \eqref{eq:f-H-0}. 
It was then shown by the first author \cite{Ro04} that $H_\mu$
is in fact a property of the phase boundary $\partial^*\{u=1\}$, see
Appendix \ref{app:defs}. This is crucial in order to apply the
(stationary) convergence 
result \cite{Sch01} to evolution problems \cite{Ro04,Ro05}.

The higher-multiplicity problem is even more challenging in the context
of the sharp interface limit of diffuse approximations,
due to the singular nature of this limit process. 
The three main ingredients of our proof are first an earlier result
of the second author \cite{To05} on the convergence of certain phase
field equations with chemical potential. This ensures rectifiability,
existence of a weak mean curvature with appropriate regularity and the
multiplicity-dependent Gibbs--Thomson relation \eqref{eq:GT-mult}. The
second ingredient is the fine local analysis of Sch\"atzle \cite{Sch04}
on rectifiable measures with sufficiently regular weak mean
curvature. The third important
argument is a comparison principle for the phase fields
$u_\eps$ and diffuse approximations of suitably constructed comparison
graphs.
%
%==================================================
%
\subsection{Organization of the paper}
\label{subsec:orga}
In the next section we will precisely formulate our assumptions and
introduce some notations. Section \ref{sec:claim} states our main
result.
A localization step in Section \ref{sec:localize}
prepares a contradiction argument that we will use in Section
\ref{sec:proof-claim} to prove our main Theorem \ref{the:claim}. We first
assume that a certain comparison
principle, which is  given in Proposition \ref{prop:comp}, holds. Section
\ref{sec:proof-comp-psi} is then devoted to the proof of this
Proposition. Finally we give in Section \ref{sec:applications} two
applications of Theorem \ref{the:main-intro} and we recall in the appendix
the definition of a generalized mean curvature for phase
boundaries that we will use.

%==================================================
%
\subsection*{Acknowledgment}
The research of Y.~Tonegawa was partially funded by the
\emph{Grant-in-aid for scientific research (B) No. 17340041}.
M.~R\"{o}ger thanks the Department of Mathematics,
Hokkaido University Sapporo for their hospitality during his visit in
August 2006. Y. Tonegawa thanks the 
Max Planck Institute for
Mathematics in the Sciences, Leipzig for their hospitality during his
visit in February 2007.

%==================================================
% assumptions
%==================================================
\section{Notations and assumptions}\label{sec:nota}
We state first all assumptions and
definitions, including those already appeared in the
introduction.
\begin{assumption}\label{Ass:1}
Consider a bounded domain $\Omega\subset\Rn$ with Lipschitz-boundary and
the standard double-well potential $W$ given by
\begin{gather*}
  W(r)\,:=\, \frac{1}{4}\big(1-r^2)^2.
\end{gather*}
We define an energy functional $E_\eps$ on $W^{1,2}(\Omega)$,
\begin{gather}
  E_\eps(u)\,:=\, \int_\Omega \Big(\frac{\eps}{2}|\nabla u|^2
  +\frac{1}{\eps}W(u)\Big)\,dx. \label{eq:def-E}
\end{gather}
Suppose  $p>n$ and let sequences $(u_\eps)_{\eps>0}\subset
W^{3,p}(\Omega)$ and $(f_\eps)_{\eps>0}\subset W^{1,p}(\Omega)$ be given
such that 
\begin{gather}
   E_\eps(u_\eps)\,\leq\, \Lambda \quad\text{ for all }\eps>0,
   \label{eq:bound-Lambda}\\
  -\eps \Delta u_\eps + \frac{1}{\eps}W^\prime(u_\eps)\,=\,
   f_\eps\quad\text{ in }\Omega. \label{eq:GT-eps}
%   u_\eps\cdot\nu_\Omega\,&=\,0\quad&&\text{ on
%   }\partial\Omega. \label{eq:neumann} 
\end{gather}
Assume further that
\begin{gather}
  u_\eps \,\to\, u\quad\text{ in } L^1(\Omega), \qquad
  u\,\in\,  BV(\Omega,\{-1,1\}), \label{ass:u}\\
  f_\eps \,\to\, f\quad\text{ weakly in }W^{1,p}(\Omega).
  \label{ass:f}
\end{gather}
\end{assumption}
We may generalize $W$ to be any $C^3$-function with two non-degenerate
minima and one local maximum, so that the results in \cite{To05} apply.

We next associate \emph{diffuse surface-area measures} and appropriate
varifolds to the functions $u_\eps$.
\begin{definition}
\label{def:vari-eps}
For $u_\eps$ we define Radon-measures $\mu_\eps$ on $\Omega$,
\begin{gather}
  \mu_\eps\,:=\, \Big(\frac{\eps}{2}|\nabla u_\eps|^2 +
  \frac{1}{\eps}W(u_\eps)\Big) \LL^n, \label{eq:def-mu-eps}
\end{gather}
and $(n-1)$-varifolds $V_\eps$ on $G^{n-1}(\Omega)$,
\begin{gather}
  V_\eps(\zeta)\,=\, \int_{G^{n-1}(\Omega)} \zeta(x,S)\,dV_\eps(x,S)\,:=\,
    \int_\Omega \zeta(x,\nu_\eps^\perp(x))\,d\mu_\eps(x)
\end{gather}
for all $\zeta\in C^0_c(G^{n-1}(\Omega))$, where
\begin{gather*}
  \nu_\eps\,:=\,\frac{\nabla u_\eps}{|\nabla u_\eps|}\text{ if }\nabla
  u_\eps\neq 0,\qquad \nu_\eps\,=\, (1,0,\dots,0)^T \text{ otherwise.}
\end{gather*}
\end{definition}
\begin{remark}
By the Sobolev embedding Theorem, \eqref{ass:f}, and $p>n$ it follows
that
\begin{gather}
  f_\eps, f\,\in\, C^{0,\alpha}(\overline{\Omega}),
  \label{eq:hoelder-f}\\
  f_\eps\,\to\, f\quad\text{ in }C^{0,\beta}(\overline{\Omega})
  \label{eq:conv-f-hoelder} 
\end{gather}
for $\alpha:=1-n/p$, all $0\leq\beta<\alpha$, and a subsequence $\eps\to
0$. 

Moreover, by \eqref{eq:bound-Lambda} there exists a subsequences
$\eps\to 0$ and a Radon-measure $\mu$ on $\Omega$ such that
\begin{gather}
  \mu_\eps\,\to\, \mu\quad\text{ as Radon-measures on }
  \Omega. \label{eq:conv-mu}
\end{gather}
Here and in the following we often do not relabel subsequences. In
particular we assume from now on that \eqref{eq:conv-f-hoelder},
\eqref{eq:conv-mu} hold for the whole sequence $\eps\to 0$.
\end{remark}

Finally we
define the mean-curvature operator $\MK$ for
graphs: for $p\in\R^{n-1}$, $X\in \Sym(n-1)$ set
\begin{gather*}
  \MK(p,X)\,:=\, (1+|p|^2)^{-\frac{3}{2}}\Big( \trace X + |p|^2\big(Id
  -\frac{p}{|p|}\otimes\frac{p}{|p|}\big):X\Big). 
\end{gather*}
%===============================================
% claim
%===============================================
\section{Statement of results}\label{sec:claim}
The first conclusion we draw is a direct consequence of previous results
of the second author.
\begin{theorem}[\cite{To05}]
\label{the:To05}
Let Assumption \ref{Ass:1} hold, let $\mu$ satisfy \eqref{eq:conv-mu},
and set
\begin{gather}\label{eq:def-sigma} 
  \sigma \,:=\, \int_{-1}^1 \sqrt{\frac{1}{2}W(s)}\,ds.
\end{gather}
Then $(2\sigma)^{-1}\mu$ is $(n-1)$-integer-rectifiable,
\begin{gather*}
  \theta^{(n-1)}(\mu,\cdot)\,=\, N(\cdot)2\sigma,\quad N \text{
 integer-valued}.
\end{gather*}
Moreover $\mu$ has weak mean curvature ${H}_\mu\in\ L^\infty(\mu)$,
and
\begin{gather}\label{eq:gt-H}
  {H}_\mu\,=\, \frac{f}{N(\cdot)\sigma}\nu
\end{gather}
holds $\mu$-almost everywhere, where $\nu=\frac{\nabla u}{|\nabla u|}$
on $\partial^*\{u=1\}$ and $\nu=0$ elsewhere.
\end{theorem}
\begin{proof}
See \cite{To05}.
\end{proof}
Our main results are summarized in the following theorem.
\begin{theorem}\label{the:claim}
Let Assumption \ref{Ass:1} hold, let $\mu$ satisfy
\eqref{eq:conv-mu}, and let $N$ denote the multiplicity function of
$\mu$ as in Theorem \ref{the:To05}.
Then
\begin{enumerate}[leftmargin=*]
\item
$\mu$-almost everywhere in $\{N(\cdot)\geq 3\text{ odd}\}$
\begin{gather}\label{eq:claim}
  H_\mu \,=\,0,\quad f\,=\, 0.
\end{gather}
\item
$\mu$-almost everywhere in $\{N(\cdot)\geq 2\text{ even}\}$
\begin{gather}
  H_\mu \,=\,0,\qquad 
  f\,\leq\, 0 \text{ in }\{u=1\}\cap\spt(\mu),\quad
    f\,\geq\, 0 \text{ in }\{u=-1\}\cap\spt(\mu).
  \label{eq:claim-even}
\end{gather}
\item
the function $H:\partial^*\{u=1\}\,\to\,\R$ defined as
\begin{gather}
  H\,:=\, H_\mu\cdot \frac{\nabla u}{|\nabla
  u|}
\end{gather}
is the generalized mean curvature of $\partial^*\{u=1\}$ in the sense of
Definition \ref{def-h}.
\item
Finally
\begin{gather}
  \sigma H \,=\, f \label{eq:GT}
\end{gather}
holds $\Ha^{n-1}$-almost everywhere on $\partial^*\{u=1\}\cap\Omega$.
\end{enumerate}
\end{theorem}
%========================================================
% prev
%========================================================
\section{Localization}
\label{sec:localize}
In this section we show that we can restrict ourselves to a
`generic' local situation, where the support of $\mu$ is well
described in terms of graphs. We then apply a result of
Sch\"atzle \cite{Sch04} that gives a fine description of the varifold
$\mu$.
\begin{lemma}
\label{lem:red1}
It is sufficient to prove \eqref{eq:claim}, \eqref{eq:claim-even} for
$\mu$-almost all 
\emph{generic points}, that are those points $x_0\in \Omega$ satisfying
\begin{align}
  & T_{x_0}\mu \text{ exists}, \label{ass:Tx0}\\
  & N_0\,:=\,N(x_0)\,\in\N, \label{ass:x0-int}\\
  & \theta\big(\mu,\{N(\cdot)=N_0\},x_0\big)\,=\,
  1. \label{ass:x0-dens}
\end{align}
In addition we may assume without loss of generality that
\begin{gather}
  x_0\,\in\,\spt(\mu)\cap \partial^*\{u=1\}\quad\text{ if }N(x_0)\text{
  is odd,} \label{ass:x0-spt}\\
  \lim_{\varrho\to 0} \Big(\sup \big\{\varrho^{-1}\dist(x,T_{x_0}\mu):
  x\in\spt(\mu)\cap B_{\varrho}^n(x_0)\big\}\Big)\,=\, 0. \label{ass:tan-sim}
\end{gather}
\end{lemma}
\begin{proof}
We show that \eqref{ass:Tx0}-\eqref{ass:tan-sim} hold
$\mu$-almost everywhere in $\{N(\cdot)\geq 2\}$.

Since $\mu$ is integer-rectifiable \eqref{ass:Tx0},
\eqref{ass:x0-int} are satisfied $\mu$-almost everywhere. Since the set 
$\{\theta^{n-1}(\mu,\cdot)=N_0\}$ is $\mu$-measurable
\eqref{ass:x0-dens} holds for $\mu$-almost all points.
By \cite[Theorem 1]{To05} $\Ha^{n-1}$-almost all
$x\in \spt(\mu)$ with odd density $N(x)$ belong to
$\partial^*\{u=1\}$. Finally \eqref{ass:tan-sim} follows from
\cite[Lemma 17.11]{Si}. 
\end{proof}
We fix $x_0$ such that \eqref{ass:Tx0}-\eqref{ass:tan-sim}
hold. After applying a suitable translation and 
rotation we may assume that $x_0=0$ and
\begin{gather}
  T_0\mu\,=\,\R^{n-1}\times\{0\}.\label{ass:Tx0-hori}
\end{gather}
To apply a contradiction argument we assume \eqref{eq:claim},
\eqref{eq:claim-even} to be false.
\begin{assumption}
\label{Ass:false}
Suppose that
\begin{gather}
  N_0\,\geq\, 2\text{ and}\quad   f(0)\,\neq\, 0.\label{ass:x0-odd}
\end{gather}
\end{assumption}
By \eqref{eq:hoelder-f}, \eqref{ass:tan-sim} there exists
$\varrho_0>0$ such that 
\begin{gather}
 |f-f(0)|\,\leq\, {\frac{1}{27}} |f(0)|\quad\text{ on
  }\overline{B_{\varrho_0}(0)}\times
  [-5\varrho_0,5\varrho_0], \label{ass:fsmall} \\
 \spt(\mu)\cap \big(B_{\varrho_0}(0)\times (-5\varrho_0,5\varrho_0)\big)\,
 \subset\, B_{\varrho_0}(0)\times
 (-\varrho_0,\varrho_0). \label{ass:spt}
\end{gather}
We distinguish four cases depending on whether $f(0)<0$ or $f(0)>0$ and
whether $u=1$ or $u=-1$ in the region `above' $\spt(\mu)$.
In the following we consider the case that
\begin{align}
 f(0)\,&>\, 0, \label{ass:f0}\\
 u(y,t)\,&=\, 1\quad\text{ for }y\in B_{\varrho_0}(0), t>\varrho_0.
 \label{ass:val-u+}
\end{align}
That implies that
\begin{gather}
  u(y,t)\,=\, \begin{cases}
  -1\quad&\text{ if } N_0\text{ is odd},\\
   1 &\text{ if }N_0\text{ is even}
 \end{cases}
 \qquad\text{ for all }y\in B_{\varrho_0}(0), t<-\varrho_0. \label{ass:val-u-}
\end{gather}
The other cases can be treated analogously or follow from a symmetry
argument. 
By \eqref{eq:conv-f-hoelder} we obtain that
\begin{gather}
 |f_\eps-f|\,\leq\, {\frac{1}{27}} f(0)\quad\text{ on
  }\overline{B_{\varrho_0}(0)}\times
  [-5\varrho_0,5\varrho_0] \label{ass:f-diff-small}
\end{gather}
for all $\eps>0$ sufficiently small.

In the next step we apply a result of Sch\"atzle \cite{Sch04} on the
local structure of the measure $\mu$. First we need some definitions.
\begin{definition}
\label{def:height}
We define the upper and lower height-functions $\varphi_+,\varphi_- :
B_{\varrho_0}^{n-1}(0)\to [-\infty,\infty]$,
\begin{align}
 \varphi_+(y)\,&:=\, \sup\big\{ t\in (-5\varrho_0,5\varrho_0)\,:\, (y,t)
  \in \spt(\mu)\cap 
  \big(B_{\varrho_0}^{n-1}(0)\times\R\big)\big\}, \label{eq:def-varphi+}\\
  \varphi_-(y)\,&:=\, \inf\big\{ t\in (-5\varrho_0,5\varrho_0)\,:\,
  (y,t) \in \spt(\mu)\cap 
  \big(B_{\varrho_0}^{n-1}(0)\times\R\big)\big\}, \label{eq:def-varphi-}
\end{align}
with the convention that the supremum over an empty set is $-\infty$ and
the infimum over an empty set is $+\infty$.

Moreover we set
\begin{gather}
  \fplus(y)  \,:=\,
  \begin{cases}
  \frac{f(y,\varphi_+(y))}{N(y,\varphi_+(y))\sigma}\quad&\text{ if }
   N(y,\varphi_+(y))\text{ is odd},\\ 
  0 \quad&\text{ if }
   N(y,\varphi_+(y))\text{ is even},
 \end{cases}
 \label{eq:def-F+}
\end{gather}
and
\begin{gather}
  \fmin(y)  \,:=\,
  \begin{cases}
  \frac{f(y,\varphi_-(y))}{N(y,\varphi_-(y))\sigma}\quad&\text{ if }
   N(y,\varphi_-(y))\text{ is odd},\\ 
  0 \quad&\text{ if }
   N(y,\varphi_-(y))\text{ is even}.
 \end{cases}
 \label{eq:def-F-}
\end{gather}
\end{definition}
\begin{proposition}
\label{prop:Sch04}
The upper height-functions $\varphi_\pm$ are twice approximately
differentiable $\LL^{n-1}$-almost 
everywhere in $\{\varphi_{\pm}\in\R\}$ with
\begin{gather}\label{eq:gt-phi}
  -\MK(\nabla\varphi_+,D^2\varphi_+) \,=\, -F_+,\\
  -\MK(\nabla\varphi_-,D^2\varphi_-) \,=\, -(-1)^{1+N_0}F_-,
  \label{eq:gt-phi-} 
\end{gather}
Moreover, for all $\frac{n-1}{2}<s<\infty$ the upper height-function
$\varphi_+$ is an $W^{2,s}$-viscosity subsolution of
\begin{gather}\label{eq:gt-odd-visc+}
  -\MK(\nabla\varphi_+,D^2\varphi_+) \,\leq\, -F_+
\end{gather}
and the lower height-function
$\varphi_-$ is a $W^{2,s}$-viscosity supersolution of
\begin{gather}\label{eq:gt-odd-visc-}
  -\MK(\nabla\varphi_-,D^2\varphi_-) \,\geq\, - (-1)^{1+N_0}F_-
\end{gather}
\end{proposition}
\begin{proof}
Since $H_\mu\in L^\infty(\mu)$ we obtain from \cite[Theorem 6.1]{Sch04} that
\begin{gather}\label{eq:gt-phi-sch}
  {H}_{\mu}(\cdot,\varphi_+(\cdot))\,=\, \nabla\cdot
  \Big(\frac{\nabla \varphi_+}{\sqrt{1+|\nabla \varphi_+|^2}}\Big)
  \frac{1}{\sqrt{1+|\nabla \varphi_+|^2}}
  \begin{pmatrix}
    -\nabla \varphi_+\\ 1
  \end{pmatrix}
\end{gather}
$\LL^{n-1}$-almost everywhere in $\{\varphi_+\in\R\}$ and that
$\varphi_+$ is for all $\frac{n-1}{2}<s<\infty$ a $W^{2,s}$-viscosity subsolution of
\begin{gather}\label{eq:gt-visc}
  -\nabla\cdot
  \Big(\frac{\nabla \varphi_+}{\sqrt{1+|\nabla \varphi_+|^2}}\Big)
  \,\leq\,-{H}_{\mu}(\cdot,\varphi_+(\cdot))\frac{1}{\sqrt{1+|\nabla \varphi_+|^2}}
  \begin{pmatrix}
    -\nabla \varphi_+\\ 1
  \end{pmatrix}
  .
\end{gather}
From \eqref{eq:gt-H}, \eqref{eq:gt-phi-sch} we deduce that
\begin{gather}
  \frac{f(\cdot,\varphi_+(\cdot))}{N(\cdot,\varphi_+(\cdot))\sigma}
  \nu(\cdot,\varphi_+(\cdot))    
  \,=\, \MK(\nabla\varphi_+,D^2\varphi_+)\frac{1}{\sqrt{1+|\nabla
  \varphi_+|^2}} 
  \begin{pmatrix}
    -\nabla \varphi_+\\ 1
  \end{pmatrix}
  . \label{eq:gt-H-recap}
\end{gather}
Next we observe that \eqref{ass:val-u+} implies
\begin{gather}
  \nu(y,\varphi_+(y))\,=\, \frac{\nabla u}{|\nabla u|}(y,\varphi_+(y))\,=\,
  \big(1+|\nabla\varphi_+|^2\big)^{-\frac{1}{2}}
  (-\nabla\varphi_+(y),1)^T \label{eq:normal-varphi+}
\end{gather}
for $\LL^{n-1}$ almost all $y\in \{\varphi_{+}\in\R\}$ such that
$(y,\varphi_+(y))\in\partial^*\{u=1\}$, and $\nu=0$ otherwise. Since up
to a $\mu$-nullset $\partial^*\{u=1\}=\{N\geq 1\text{ odd}\}$ we
obtain from \eqref{eq:def-F+} and \eqref{eq:gt-H-recap},
\eqref{eq:normal-varphi+} that \eqref{eq:gt-phi}
holds. \eqref{eq:gt-odd-visc+} follows by the same arguments.
To obtain \eqref{eq:gt-phi-}, \eqref{eq:gt-odd-visc-} we observe that
\begin{gather}
  \nu(y,\varphi_-(y))\,=\, \frac{\nabla u}{|\nabla
  u|}(y,\varphi_-(y))\,=\, (-1)^{1+N_0} 
  \big(1+|\nabla\varphi_-|^2\big)^{-\frac{1}{2}}
  (-\nabla\varphi_-(y),1)^T \label{eq:normal-varphi-}
\end{gather}
and we proceed as above.
\end{proof}
We choose below a `good point' for which we derive a contradiction to
Assumption \ref{Ass:false}.
Before, we need another definition.
\begin{definition}
\label{def:taylor}
We say that a function $\psi$ has a second-order Taylor expansion
at a point $y_1\in\R^{n-1}$ if there exist $p\in\R^{n-1}, X\in
\Sym(n-1)$ such that
\begin{gather}
  \psi(y)-\psi(y_1) -p\cdot(y-y_1)
  -\frac{1}{2}(y-y_1)\cdot X(y-y_1) \,=\,
  o(|y-y_1|^2). \label{eq:taylor} 
\end{gather}
We then set $ \nabla\psi(y_1):=p$, $D^2\psi(y_1):=X$.
\end{definition}
\begin{lemma}
\label{lem:red2}
There exists a point $y_1\in B_{\varrho_0}(0)$ such
that 
\begin{align}
  &\varphi_+(y_1)\,=\,\varphi_-(y_1), \label{eq:phi-y1}\\
  &x_1\,:=\, (y_1,\varphi_\pm(y_1)) \text{ is a generic point},
  \label{eq:generic-y1}\\
  &\theta^{n-1}(\mu,x_1)\,=\, N_0, \label{eq:dens-y1}\\
  &\varphi_\pm\text{ have a second-order Taylor expansion at }y_1,
  \label{eq:phi-taylor}\\
  & -\MK(\nabla\varphi_\pm(y_1),D^2\varphi_\pm(y_1))\,=\,
  -F_\pm(y_1). \label{eq:phi-MK} 
\end{align}
\end{lemma}
\begin{proof}
Since the weak mean curvature $H_\mu$ belongs to $L^\infty(\mu)$ by
\eqref{eq:gt-H} and since $0\in\Rn$ is a generic point we can apply
\cite[Lemma 3.4]{Sch01} (see also Step 3 in the proof of \cite[Lemma
3.2]{Ro05}) and obtain that the set
\begin{gather*}
  \Sigma_0\,:=\, \Big\{x=(y,\varphi_\pm(y)) : y\in
  B^{n-1}_{\varrho_0}(0)\cap \{\varphi_+=\varphi_-\},\qquad\qquad\\
  \,\qquad\qquad x\in\spt(\mu)\text{ is generic, } N(x)\,=\,N_0\Big\} \label{eq:S0}
\end{gather*}
has full density with respect to $\mu$ in $0\in\Rn$. This property was
essentially deduced from a tilted version of Brakke's Lipschitz 
Approximation Theorem. The curvature bound ensures a strong control
on the approximations, see \cite{Sch01,Ro04} for the details.
From the Coarea Formula we then deduce that
\eqref{eq:phi-y1}-\eqref{eq:dens-y1} holds in a set with full
$\LL^{n-1}$-density in $0\in\R^{n-1}$. 
Finally, \eqref{eq:phi-taylor}, \eqref{eq:phi-MK} are satisfied
$\LL^{n-1}$ almost everywhere in $\{\varphi_+=\varphi_-\}$ by
\cite[Proposition 4.1]{Sch04} and \eqref{eq:gt-phi},
\eqref{eq:gt-phi-}. 
\end{proof}
%============================================
%  proof claim
%============================================
\section{Proof of Theorem \ref{the:claim}}
\label{sec:proof-claim}
We fix $y_1\in B_{\varrho_0}(0)$ such that
\eqref{eq:phi-y1}-\eqref{eq:phi-MK} hold and consider the second
order Taylor approximation of $\varphi_\pm$ at $y_1$,
\begin{gather}
  P_1(y)\,:=\, \varphi_\pm(y_1) +\nabla\varphi_\pm(y_1)\cdot (y-y_1)
  +\frac{1}{2}(y-y_1)\cdot D^2\varphi_\pm(y_1)(y-y_1). \label{eq:def-P1}
\end{gather}
From \eqref{eq:phi-taylor} we then deduce that
\begin{gather}
  |P_1(y) - \varphi_\pm(y)|\,=\, o(|y-y_1|^2). \label{eq:taylor-y1}
\end{gather}
\begin{lemma}
\label{lem:ex-psi}
For all $\smallbdry>0$ there is $\varrho_1>0$ such that
\begin{gather}
  B_{\varrho_1}(y_1)\subset B_{\varrho_0}(0), \label{eq:B-rho1}\\
  \big| P_1(y)- \varphi_+(y)\big| \,<\,\smallbdry |y-y_1|^2 \quad\text{ on }
  B_{\varrho_1}(y_1) \label{eq:comp-psi-bdry}
\end{gather}
and such that for any
$0<\varrho<\varrho_1$ there exists a unique solution $\psi\in
C^\infty(\overline{B_{\varrho}(y_1)})$ of
\begin{alignat}2
  -\MK(\nabla\psi,D^2\psi)\,&=\,
  -\frac{2f(y_1)}{3\sigma}\quad&&\text{ in }B_{\varrho}(y_1),
  \label{eq:psi}\\
   \psi(y) \,&=\, P_1(y) +\smallbdry\varrho^2\quad&&\text{ on }\partial
  B_{\varrho}(y_1). \label{eq:psi-bdry} 
\end{alignat}
\end{lemma}
\begin{proof}
Since $y_1\in B_{\varrho_0}(0)$ and by \eqref{eq:taylor-y1} for any
$\varrho_1>0$ sufficiently small the properties \eqref{eq:B-rho1},
\eqref{eq:comp-psi-bdry} hold. Since $P_1$ is smooth and since the
right-hand side of equation \eqref{eq:psi} is constant we deduce from
\cite[Theorem 16.9]{GT} that for all
\begin{gather*}
  \varrho\,\leq\,\varrho_1,\quad
  \varrho_1=\varrho_1(n,\sigma,f(y_1))
\end{gather*}
a unique solution $\psi\in C^{2,\gamma}(\overline{B_{\varrho}(y_1)})$,
$0<\gamma<1$, of 
\eqref{eq:psi}, \eqref{eq:psi-bdry} exists. The higher
regularity of $\psi$ follows from standard elliptic theory and the
smoothness of the data in \eqref{eq:psi}, \eqref{eq:psi-bdry}.
\end{proof}
The next Proposition is the heart of the contradiction argument. It
relies on the fact that the approximations $u_\eps$ behave as if the
curvature of the limit interface is given by $f/\sigma$ rather than
by $f/(N\sigma)$. 
\begin{proposition}
\label{prop:comp-psi}
Let $\psi$ be as in Lemma \ref{lem:ex-psi}. Then
\begin{gather}
  \psi\,\geq\, \varphi_+ \quad\text{ in }B_{\varrho}(y_1).
  \label{eq:comp-psi}
\end{gather}
\end{proposition}
The proof of this Proposition uses a comparison between
$u_\eps$ and approximations $v_\eps$ of $2\Chi_E-1$, where $E$ is the
region above the graph of $\psi$. We postpone this proof to section
\ref{sec:proof-comp-psi} and continue the proof of 
Theorem \ref{the:claim}.
\begin{lemma}
\label{lem:eta}
For all $0<\smalleta<1$ there exists $\tilde{\varrho}_1>0$ such that for
all $0<\varrho<\tilde{\varrho}_1$ the function 
\begin{gather}
  \eta(y)\,:=\, P_1(y) - 2\smalleta(\varrho^2-|y-y_1|^2)
  \label{eq:def-eta}
\end{gather}
satisfies for all $y\in B_{\varrho}(y_1)$
\begin{gather}
  -\MK(\nabla\eta,D^2\eta)\,\geq\,
  -F_\pm(y_1) -5(n-1)\smalleta.
  \label{eq:eta}
\end{gather}
\end{lemma}
\begin{proof}
We compute that for $y\in B_\varrho(y_1)$
\begin{align*}
  |\nabla\eta(y) - \nabla\varphi_\pm(y_1)|\,&=\,
   |D^2\varphi_\pm(y_1)(y-y_1)+4\smalleta(y-y_1)|\\
   &\leq\,
   \big(|D^2\varphi_\pm(y_1)|+4\big)\varrho,\\
   |D^2\eta(y)|\,&\leq\, |D^2\varphi_\pm(y_1)|+4.
\end{align*}
Hence, we can choose
$\tilde{\varrho}_1=
\tilde{\varrho}_1(\nabla\varphi_\pm(y_1),D^2\varphi_\pm(y_1),\smallbdry)$   
such that
\begin{align*}
  |\MK(\nabla\eta,D^2\eta) - \MK(\varphi_\pm(y_1),D^2\eta)|\,\leq\,
   (n-1)\smalleta\quad\text{ in }B_\varrho(y_1)
\end{align*}
for all $0<\varrho<\tilde{\varrho}_1$. This implies that 
\begin{align*}
  -\MK(\nabla\eta,D^2\eta)\,&\geq\,
  -\MK(\nabla\varphi_\pm(y_1),D^2\varphi_\pm(y_1))
  -4\smalleta\MK(\nabla\varphi_\pm(y_1),Id) - (n-1)\smalleta\\
  &\geq -\MK(\nabla\varphi_\pm(y_1),D^2\varphi_\pm(y_1))
  - 5(n-1)\smalleta,
\end{align*}
where we have used that $\MK(p,Id)\leq n-1$ for all $p\in\R^{n-1}$.
\end{proof}
\begin{proof}[{\bf Proof of Theorem \ref{the:claim}}]
Choose $0<\smallbdry<\frac{f(y_1)}{30(n-1)\sigma}$ and
$0<\varrho<\min(\varrho_1,\tilde{\varrho}_1)$.
Let $\psi,\eta$ be the functions constructed in Lemma \ref{lem:ex-psi} and Lemma
\ref{lem:eta}.
We then obtain from \eqref{eq:psi}, \eqref{eq:eta}, the definition of
$F_\pm$, \eqref{ass:x0-odd}, and $f(y_1)>0$ that
\begin{align}
  -\MK(\nabla\psi,D^2\psi) + \MK(\nabla\eta,D^2\eta)\,&\leq\,
  -\frac{2f(y_1)}{3\sigma} 
  +\frac{f(y_1)}{N(y_1,\varphi_\pm(y_1))\sigma} + 5(n-1)\smalleta\notag\\
  &\leq\,  -\frac{2f(y_1)}{3\sigma} 
  +\frac{f(y_1)}{2\sigma} + 5(n-1)\smalleta\notag\\
  &<\, 0.
  \label{eq:psi-eta} 
\end{align}
Since $|\nabla\psi|,|\nabla\eta|$ are uniformly bounded  the
maximum principle \cite[Theorem 10.1]{GT} implies that $\psi-\eta$ has
no interior maximum. 
In particular,
\begin{align*}
  \psi(y_1)-\varphi_\pm(y_1) +2\smalleta\varrho^2 \,&=\,
  \psi(y_1)-\eta(y_1)\\
  &\leq\, \sup_{\partial B_\varrho(y_1)}(\psi-\eta)\,=\, \smallbdry
  \varrho^2 
\end{align*}
and we deduce that
\begin{gather*}
  \psi(y_1)-\varphi_\pm(y_1) \,<\,0,
\end{gather*}
which is a contradiction to \eqref{eq:comp-psi}.

This shows by Assumption \ref{Ass:false} and \eqref{ass:f0},
\eqref{ass:val-u+} that
\begin{gather}
  f\,\leq\, 0\quad\text{ if }\, N_0\geq 2\text{ and } u=1\, \text{ in }
  \{t>\varphi_+(y)\}. \label{eq:con-1}
\end{gather}
By a symmetry argument $u\mapsto -u$ it follows that
\begin{gather}
  f\,\geq\, 0\quad\text{ if } N_0\geq 2\text{ and } u=-1 \text{ in }
  \{t>\varphi_+(y)\}. \label{eq:con-2}
\end{gather}
As we explain in Remark \ref{rem:symmetry} below we obtain also
\begin{gather}
  f\,\geq\, 0\quad\text{ if } N_0\geq 3\text{ is odd and }u=1 \text{ in }
  \{t>\varphi_+(y)\} \label{eq:con-3}
\end{gather}
and, again by symmetry, that
\begin{gather}
  f\,\leq\, 0\quad\text{ if } N_0\geq 3\text{ is odd and }u=-1 \text{ in }
  \{t>\varphi_+(y)\}. \label{eq:con-4}
\end{gather}
Putting together \eqref{eq:con-1}-\eqref{eq:con-4} and using
\eqref{eq:gt-H} we deduce the conclusion (1), (2) of Theorem
\ref{the:claim}. The conclusion (3) follows from Proposition
\ref{prop-def-h}, and the statement (4) is deduced from (1), (3) and
\eqref{eq:gt-H}. 
\end{proof}
\begin{remark}\label{rem:symmetry}
In the case that $f(0)<0$ and $u=1$ `above' $\varphi_+$ one
considers for suitably small $\varrho>0$ the  solution $\tilde{\psi}$ of
\begin{alignat*}2
  -\MK(\nabla\tilde{\psi},D^2\tilde{\psi})\,&=\,
  -\frac{2f(y_1)}{3\sigma}\quad&&\text{ in }B_{\varrho}(y_1),
  \\
   \tilde{\psi}(y) \,&=\, P_1(y) -\smallbdry\varrho^2\quad&&\text{ on }\partial
  B_{\varrho}(y_1)
\end{alignat*}
and the function $\tilde{\eta}$,
\begin{gather*}
  \tilde{\eta}(y)\,:=\, P_1(y) + 2\smalleta(\varrho^2-|y-y_1|^2)
\end{gather*}
To derive a contradiction the corresponding statement to Proposition
\ref{prop:comp-psi} is needed, that is 
\begin{gather*}
  \tilde{\psi}\,\leq\, \varphi_- \quad\text{ in }B_{\varrho}(y_1).
\end{gather*}
In the case that $N_0\geq 1$ this property can be proved in the same way
as we will prove Proposition \ref{prop:comp-psi}: One constructs 
smooth approximations $v_\eps$ of the function $2\Chi_{\{t>\tilde{\psi}(y)\}}-1$
and uses a comparison principle to obtain $u_\eps\leq v_\eps$.
However, these arguments do not apply if $N_0\geq 2$ is even,
since in that case $u_\eps\approx 1$ is larger than
$2\Chi_{\{t>\tilde{\psi}(y)\}}-1$ in the region `below' $\tilde{\psi}$.
\end{remark}
%===============================================
% proof of prop
%===============================================
\section{Proof of Proposition \ref{prop:comp-psi}}
\label{sec:proof-comp-psi}
Assume that \eqref{eq:comp-psi} does not hold, that is 
\begin{gather}
  \sup_{B_\varrho(y_1)} (\varphi_+ -\psi)\,>\,0. \label{eq:ass-contra}
\end{gather}
By \eqref{eq:comp-psi-bdry} and \eqref{eq:psi-bdry}
\begin{gather}
  \varphi_+\,<\,\psi\quad\text{ on }\partial
  B_\varrho(y_1). \label{eq:comparison-bdry1}
\end{gather}
Since $\psi$ is continuous and $\varphi_+$ is upper-semicontinuous there
exist $0<\varrho_3<\varrho_2<\varrho$ such that
\begin{gather}
  \sup_{B_{\varrho_3}(y_1)} (\varphi_+-\psi)\,>\, 0, \label{eq:int-max}\\
  \varphi_+\,<\,\psi\quad\text{ on } B_\varrho(y_1)\setminus
  B_{\varrho_3}(y_1). \label{eq:psi-phi-bdry} 
\end{gather}
As explained before we will use that $u_\eps$ behaves as if the
curvature of the sharp interface limit is given by $f/\sigma$, instead
of $f/(N\sigma)$. In a first step we construct functions
$v_\eps$ such that 
\begin{gather}
  -\eps \Delta v_\eps + \frac{1}{\eps}W^\prime(v_\eps)\,\leq\,
   \frac{7}{9}f(y_1) \label{eq:v-eps},\\
   v_\eps\,\to\, 2\Chi_{\{(y,t):t>\psi(y)\}}-1. \label{conv:v}
\end{gather}
In the second step we will apply a comparison principle to
$u_\eps,v_\eps$ to obtain a contradiction in the limit $\eps\to 0$.
\subsection{Construction of $v_\eps$}
The two ingredients to construct $v_\eps$ are a \emph{modified distance}
function from $\graph(\psi)$ and the optimal profile and first
order-correction of the one-dimensional minimisation problem associated
to the Cahn--Hilliard functional.
\begin{definition}
\label{def:M}
We define
\begin{gather*}
  M\,:=\,\graph\big(\psi\lfloor B_\varrho(y_1)\big)
\end{gather*}
and denote by $\ddist:=\dist(M,\cdot)$ the 
signed distance function from $M$, taken positive in the region `above'
$M$. Moreover we let $\Pi_M:\Rn\to M$ be  the orthogonal projection onto
$M$ and  $(\kappa_i)_{i=1,...,n-1}$ the principal curvatures of $M$.
Finally we define for $x\in\R^n$, $(y,t)=\Pi_M(x)$
\begin{align*}
  \tilde{\kappa}_{i}(x)\,&:=\,\kappa_i(y,t),
\end{align*}
which is well-defined in a neighborhood of $M$.
\end{definition}
\begin{remark}
\label{rem:dist-fun}
Since $\psi$ is smooth we deduce that $M$ is a smooth hypersurface and
that there exists $\delta>0$,
$\delta=\delta(\|\psi\|_{C^2(B_\varrho(y_1))},\varrho_2)$ such that
the distance function $\ddist$ is unique and smooth
in a neighborhood
\begin{gather}
  G\,=\, \big\{x \in B_{\varrho_2}(y_1)\times (-5\varrho_0,5\varrho_0) :
  |\ddist(x)|<\delta\big\}. \label{eq:def-G}
\end{gather}
Moreover
\begin{align} 
  \Delta \ddist\,=\,& \sum_{i=1}^{n-1}\tilde{\kappa}_{i}
  +O(\ddist)
  =\, \frac{2f(y_1)}{3\sigma} 
  +O(\ddist).
  \label{eq:delta-d2}
\end{align}
holds  in $G$ \cite[Lemma 14.17]{GT}.
\end{remark}
We turn to the optimal profile for the one-dimensional minimisation in
the Cahn--Hilliard energy.
\begin{remark}
\label{rem:opt-profile}
Let $\phi_0:\R\to [-1,1]$ be the \emph{optimal profile}, that is the
solution of
\begin{gather}
  -\phi_0^{\prime\prime} + W^\prime(\phi_0)\,=\, 0, \label{eq:phi}\\
  \phi_0(-\infty)\,=\, -1, \quad \phi_0(+\infty)\,=\,1,
  \label{eq:infty-phi_0}
\end{gather}
and let $\phi_1:\R\to\R$ be the first order correction (see \cite{Pao97}),
\begin{gather}
  -\phi_1^{\prime\prime} +W^{\prime\prime}(\phi_0)\phi_1\,=\, \phi_0^\prime +
  \sigma, \label{eq:phi_1}\\
  \phi_1(\pm\infty)\,=\, \frac{\sigma}{W^{\prime\prime}(\pm
    1)}. \label{eq:infty-phi_1}
\end{gather}
\end{remark}
Since the distance function $\ddist$ is smooth only in a neighborhood
of $\graph(\psi)$ we have to modify the distance function.
\begin{definition}
\label{def:mod-dist}
For $\eps>0$ we choose $\delta(\eps)>0$ such that
\begin{gather}
  \delta(\eps)\,\to 0,\, \frac{\delta(\eps)}{\eps}\,\to\,\infty\quad
  \text{ as }\eps\to 0
  \label{eq:cond-delta}
\end{gather}
and such that the conditions
\begin{gather}
  \frac{1}{\delta(\eps)}\phi_0^\prime\big(\eps^{-1}\delta(\eps)\big),\,
  \frac{1}{\eps}\phi_0^{\prime\prime}\big(\eps^{-1}\delta(\eps)\big)\,
  \to 0\quad\text{ as }\eps\to 0
  \label{eq:cond-delta-3}
\end{gather}
are satisfied.
Moreover we choose smooth functions $\beta_\eps:\R\to\R$, $\eps>0$, with
\begin{gather}
  \beta_\eps(r)\,=\,
  \begin{cases}
    r &\text{ for }\quad |r|\leq \frac{\delta(\eps)}{3}\\
    -\delta(\eps)&\text{ for }\quad r\leq -2\delta(\eps)\\
    \delta(\eps)&\text{ for }\quad r\geq 2\delta(\eps)
  \end{cases}
  \label{eq:def-beta}
\end{gather}
such that
\begin{align}
  &0\,\leq\, \beta_\eps^\prime\,\leq\, 1,\label{eq:bounds-beta}\\
  &0\,\geq\, \beta_\eps^{\prime\prime}\,\geq\,
  -\frac{3}{\delta(\eps)}. \label{eq:bounds-beta-prime}
\end{align}
We then define the \emph{modified distance functions} $d_\eps$,
\begin{gather}
  d_\eps(x)\,:=\, \beta_\eps(\ddist(x)).\label{eq:def-d}
\end{gather}
\end{definition}
\begin{remark}
\label{rem:mod-dist}
We observe that by \eqref{eq:cond-delta}
\begin{gather}
  \{|\ddist|\leq 2\delta(\eps)\}\subset G \label{eq:G}
\end{gather}
for all $\eps\leq\eps_0(\|\psi\|_{C^2(B_\varrho(y_1))})$ and deduce that
$d_\eps$ is smooth for sufficiently small $\eps>0$.
We compute that
\begin{align}
  \nabla d_\eps\,=&\, \beta_\eps^\prime(\ddist)\nabla\ddist,\label{eq:grad-d}\\
  \Delta d_\eps\,=&\, \beta_\eps^{\prime\prime}(\ddist) +
  \beta_\eps^\prime(\ddist) \Delta\ddist\notag\\
  =&\, \beta_\eps^{\prime\prime}(\ddist) +
  \beta_\eps^\prime(\ddist) \Big(\frac{2f(y_1)}{3\sigma} 
  +O(\ddist)\Big),
  \label{eq:Deltad}
\end{align}
where we have used \eqref{eq:delta-d2}. For $\eps<\eps_0(\psi)$ we
obtain that
\begin{gather*}
  \beta_\eps^\prime\circ\ddist \,=\, 0\quad\text{ in the set }\big\{|\ddist|
  > 2\delta(\eps)\}\big\}
\end{gather*}
and we deduce from \eqref{ass:fsmall}, \eqref{eq:bounds-beta},
\eqref{eq:bounds-beta-prime} that 
\begin{gather}
  |\Delta d_\eps|\,\leq\,
   C\big(1+\frac{1}{\delta(\eps)}\big). \label{eq:est-Deltad}
\end{gather}
\end{remark}
We are now ready to define $v_\eps$.
\begin{definition}
\label{def:v-eps}
Let $\eps_0=\eps_0(\|\psi\|_{C^2(B_{\varrho}(y_1))})$ be chosen such that
\eqref{eq:G} holds. We then define
$v_\eps:B_{\varrho_2}(y_1)\times\R\,\to\, \R$,
\begin{gather}
  v_\eps(x)\,:=\, \phi_0(\eps^{-1}d_\eps(x))
  +\eps\phi_1(\eps^{-1}d_\eps(x)) \frac{2}{3\sigma}f(y_1).
  \label{eq:def-v-eps}
\end{gather}
\end{definition}
\subsection{Comparison of $u_\eps, v_\eps$}
\subsubsection{Subsolution property}
We are going to show that $v_\eps$ is a suitable subsolution of a
(diffuse) constant curvature equation.

We first compute that, using \eqref{eq:phi_1},
\begin{align}
  -\eps\Delta v_\eps \,=\,&
  -\frac{1}{\eps}\phi_0^{\prime\prime}(\eps^{-1}d_\eps) |\nabla d_\eps|^2
  -\phi_0^\prime(\eps^{-1}d_\eps) \Delta d_\eps\notag\\
  &+ \frac{2f(y_1)}{3\sigma}\Big(
  -W^{\prime\prime}(\phi_0(\eps^{-1}d_\eps)) \phi_1(\eps^{-1}d_\eps) +
  \phi_0^\prime(\eps^{-1}d_\eps) + \sigma\Big)|\nabla d_\eps|^2 \notag\\
  &  -\eps\phi_1^\prime(\eps^{-1}d_\eps) \frac{2f(y_1)}{3\sigma}\Delta
  d_\eps 
  \label{eq:delta-v-eps}  
\end{align}
and
\begin{align}
  \frac{1}{\eps} W^\prime(v_\eps)\,=\,&
  \frac{1}{\eps}W^\prime(\phi_0(\eps^{-1}d_\eps)) +
  W^{\prime\prime}(\phi_0(\eps^{-1}d_\eps))
  \frac{2f(y_1)}{3\sigma}\phi_1(\eps^{-1}d_\eps) +O(\eps).
  \label{eq:W-v-eps}
\end{align}
Using \eqref{eq:est-Deltad} we
deduce from \eqref{eq:delta-v-eps}, \eqref{eq:W-v-eps} that
\begin{align}
  -\eps\Delta v_\eps +\frac{1}{\eps}W^\prime(v_\eps)
  =\,&-\frac{1}{\eps}\phi_0^{\prime\prime}\big(|\nabla
  d_\eps|^2-1\big)
  -\phi_0^\prime\big(\Delta d_\eps -\frac{2f(y_1)}{3\sigma}|\nabla
  d_\eps|^2\big)\notag\\
  &-\frac{2f(y_1)}{3\sigma}\phi_1 W^{\prime\prime}(\phi_0)(|\nabla
  d_\eps|^2-1)\notag\\
  &+\frac{2}{3}f(y_1)|\nabla d_\eps|^2 +
  o(\eps\delta(\eps)^{-1}).\label{eq:gt-v-eps} 
\end{align}
\begin{proposition}
\label{prop:v-eps}
For all $0<\eps<\eps_0$, $\eps_0=\eps_0(\psi,W)$,
\begin{gather}
  -\eps\Delta v_\eps +\frac{1}{\eps}W'(v_\eps)\,\leq\,
   \frac{7}{9}f(y_1). \label{eq:gt-v-eps-global}
\end{gather}
\end{proposition}
\begin{proof}
We check \eqref{eq:gt-v-eps-global} in the different regions.
\step{1}
In the region $\{|\ddist|\leq \frac{1}{3}\delta(\eps)\}$ holds $|\nabla
d_\eps|=1$ and we obtain from
\eqref{eq:Deltad}, \eqref{eq:gt-v-eps} that
\begin{align*}
  -\eps\Delta v_\eps +\frac{1}{\eps}W^\prime(v_\eps)
  =\,& \frac{2}{3}f(y_1)
  + O(\ddist) + o(\eps\delta(\eps)^{-1})\\
  \leq\,&  \frac{2}{3}f(y_1) +O(\delta(\eps)) + o(\eps\delta(\eps)^{-1}).
\end{align*}
Therefore \eqref{eq:gt-v-eps-global} holds for $\eps>0$
sufficiently small.
\step{2}
In $\{\ddist \geq 2\delta(\eps)\}$ we obtain
\begin{align}
  &-\eps\Delta v_\eps +\frac{1}{\eps}W^\prime(v_\eps)\notag\\
  =\,& \frac{1}{\eps}W^\prime\Big(\phi_0(\eps^{-1}\delta(\eps))
  +\eps\phi_1(\eps^{-1}\delta(\eps)) \frac{2}{3\sigma}f(y_1)
  \Big)\notag\\
  =\,& \frac{1}{\eps}W^\prime\big(\phi_0(\eps^{-1}\delta)\big)
  +W^{\prime\prime}\big(\phi_0(\eps^{-1}\delta)\big)
  \phi_1(\eps^{-1}\delta)\frac{2}{3\sigma}f(y_1) +
  O(\eps)\phi_1(\eps^{-1}\delta). 
\end{align}
From \eqref{eq:phi}, \eqref{eq:cond-delta-3} and \eqref{eq:infty-phi_1},
\eqref{eq:cond-delta} we
deduce that \eqref{eq:gt-v-eps-global} holds in  $\{\ddist \geq
2\delta(\eps)\}$ for sufficiently small $\eps>0$. By
similar calculations we  obtain \eqref{eq:gt-v-eps-global} also in
the region $\{\ddist \leq -2\delta(\eps)\}$.
\step{3}
Let us now consider the set $\{\frac{1}{3}\delta(\eps)\leq \ddist \leq
2\delta(\eps)\}$ and estimate the different terms in
\eqref{eq:gt-v-eps}. We first obtain from \eqref{eq:cond-delta-3} that
in this region  
\begin{align}
  -\frac{1}{\eps}\phi_0^{\prime\prime}(\eps^{-1}d_\eps)\big(|\nabla
   d_\eps|^2-1\big)\,\to\, 0 \quad\text{ as }\eps\to 0.
   \label{eq:gt-trans-t1}
\end{align}
Next we compute that, using \eqref{eq:grad-d}, \eqref{eq:Deltad} and
\eqref{eq:bounds-beta-prime},
\begin{align*}
  &\Big|-\phi_0^\prime(\eps^{-1}d_\eps)\big(\Delta d_\eps
  -\frac{2f(y_1)}{3\sigma}|\nabla 
  d_\eps|^2\big)\Big| \notag\\
  =\,
  &\phi_0^\prime(\eps^{-1}d_\eps)\Big|\beta_\eps^{\prime\prime}(\ddist) +
  \beta_\eps^\prime(\ddist) \Big(\frac{2f(y_1)}{3\sigma} 
  +O(\ddist)\Big) -\frac{2f(y_1)}{3\sigma}\beta_\eps^\prime(\ddist)^2
  \Big|\notag\\
  \leq\,
    &\phi_0^\prime((3\eps)^{-1}\delta(\eps))\Big(\frac{3}{\delta(\eps)} +
    \frac{4f(y_1)}{3\sigma} +   O(\delta(\eps))\Big). 
\end{align*}
Hence, by \eqref{eq:cond-delta-3}
\begin{gather}
  -\phi_0^\prime(\eps^{-1}d_\eps)\big(\Delta d_\eps
  -\frac{2f(y_1)}{3\sigma}|\nabla 
  d_\eps|^2\big)\,\to\, 0\, \text{ as }\eps\to 0. \label{eq:gt-trans-t2}
\end{gather}
Finally we observe that in $\{\frac{1}{3}\delta(\eps)\leq \ddist \leq
 2 \delta(\eps)\}$
\begin{gather}
  a_\eps\,:=\,  \frac{\phi_1(\eps^{-1}d_\eps)
  W^{\prime\prime}(\phi_0(\eps^{-1}d_\eps))}{\sigma}\,=\, 1+o(1)
\end{gather}
by \eqref{eq:infty-phi_1}, \eqref{eq:cond-delta} and we deduce for the
last line in \eqref{eq:gt-v-eps} that
\begin{gather}
  \frac{2}{3}f(y_1)\Big((-a_\eps +1)|\nabla d_\eps|^2 + a_\eps\Big)\,= \,\frac{2}{3}f(y_1)+
  o(1).\label{eq:gt-trans-t3} 
\end{gather}
We obtain from \eqref{eq:gt-v-eps} and
\eqref{eq:gt-trans-t1}, \eqref{eq:gt-trans-t2}, \eqref{eq:gt-trans-t3}
that \eqref{eq:gt-v-eps-global} holds in $\{\frac{1}{3}\delta(\eps)\leq
\ddist \leq 2 \delta(\eps)\}$ for all $\eps>0$ sufficiently small.
By similar considerations we prove \eqref{eq:gt-v-eps-global} also in
$\{-\delta(\eps)\leq \ddist \leq -\frac{1}{3}\delta(\eps)\}$.
\end{proof}
%
%====================      comparison: bulk
%
\subsubsection{Comparison in the bulk regions}
\begin{lemma}\label{lem:To}
As $\eps\to 0$
\begin{gather}
  u_\eps\,\to\, 1 \quad\text{ uniformly on each compact subset of
  }\{t>\varphi_+(y)\}, \label{eq:To1}\\
  u_\eps\,\to\, (-1)^{N_0} \quad\text{ uniformly on each compact subset of
  }\{t<\varphi_-(y)\}. \label{eq:To-1}
\end{gather}
and
\begin{gather}
  v_\eps\,\to\, 1 \quad\text{ uniformly on each compact subset of
  }\{t>\psi(y)\}, \label{eq:To-v1}\\
  v_\eps\,\to\, -1 \quad\text{ uniformly on each compact subset of
  }\{t<\psi(y)\}. \label{eq:To-v-1}
\end{gather}
\end{lemma}
\begin{proof}
By \cite[Proposition 4.2]{To05} either $u_\eps \to 1$ or $u_\eps\to -1$
uniformly on each compact subset of $\Omega\setminus \spt\mu$. By
\eqref{ass:val-u+}, \eqref{ass:val-u-} we conclude \eqref{eq:To1}, \eqref{eq:To-1}.  The
construction of $v_\eps$ yields \eqref{eq:To-v1},
\eqref{eq:To-v-1}. 
\end{proof}
For $\eps>0$ such that \eqref{ass:fsmall}, \eqref{ass:f-diff-small}
holds we deduce that
\begin{gather}
  f_\eps\,>\, \frac{8}{9}f(y_1)\quad\text{ in }
  B_{\varrho_0}(0)\times(-5\varrho_0,5\varrho_0). \label{eq:rhs-u} 
\end{gather}
In fact, in this region we compute that
\begin{gather*}
  \frac{f(y_1)}{f_\eps}\,\leq\,\frac{|f(y_1)-f(0)|+f(0)}{f(0)- |f-f(0)|-|f_\eps-f|}
  \,\leq\, \frac{28}{25}\,<\,\frac{9}{8}.
\end{gather*} 
\begin{lemma}\label{lem:bdry-1}
Choose $c_0>0$ such that
\begin{gather}
  W^{\prime\prime}(r)\,\geq\, \frac{W''(1)}{2}>0\quad\text{ for all }|r|\geq
  1-c_0 \label{eq:def-c0} 
\end{gather}
and set $\Omega':= B_{\varrho_0}(0)\times(-5\varrho_0,5\varrho_0)$.

Let $\lambda_{\eps,+}, \lambda_{\eps,-}$ be the positive and negative
solution, respectively, of  
\begin{gather}
  \frac{1}{\eps}W'(\lambda_{\eps,\pm})\,=\,
  \frac{8}{9}f(y_1). \label{eq:def-lambda-eps} 
\end{gather}
Then there exists for any bounded domain $U\subset\subset \Omega'$ and
any $k\in\N$ a constant $C_k=C_k(n,U,\Omega')$ such that for all
$\eps<\eps_0(k,c_0)$
\begin{gather}
  u_\eps\,\geq\,
  \lambda_{\eps,-} - C_k\eps^k \quad\text{ in }U. \label{eq:asym-u-eps-2}
\end{gather}
Assume \eqref{eq:def-c0}, \eqref{eq:def-lambda-eps} and in addition that
there exists $\Omega''\subset\Omega'$ such that 
\begin{gather}
  u_\eps\,\geq\, 1-c_0\quad\text{ in }\Omega''. \label{eq:u-out}
\end{gather}
Then  there exists for any bounded domain $U\subset\subset \Omega''$ and
any $k\in\N$ a constant $C_k=C_k(n,U,\Omega'')$ such that
\begin{gather}
  u_\eps\,\geq\,
  \lambda_{\eps,+} - C_k\eps^k \quad\text{ in }U. \label{eq:asym-u-eps}
\end{gather}
\end{lemma}
\begin{proof}
We first show the second conclusion.
We deduce from \eqref{eq:rhs-u}, \eqref{eq:def-lambda-eps} that
\begin{gather}
  -\eps\Delta(u_\eps-\lambda_{\eps,+})
   +\frac{1}{\eps}\Big(W'(u_\eps)-W'(\lambda_{\eps,+})\Big) \,\geq\, 0.
   \label{eq:lambda-diff}
\end{gather}
Consider first bounded domains $\Omega_1,\Omega_2$ such that
\begin{gather*}
  U\,\subset\,\Omega_1 \,\subset\,\subset\, \Omega_2\,\subset \Omega''.
\end{gather*}
Chose a cut-off function $\phi\in C^\infty_c(\Omega'')$ such that $0\leq
\phi\leq 1$ and 
\begin{gather*}
  \phi\,=\, 1 \,\text{ in } \Omega_1,\qquad
  \phi\,=\, 0 \,\text{ in } \Omega''\setminus \Omega_2.
\end{gather*}
Next we define
\begin{gather*}
  (u_\eps-\lambda_{\eps,+})_- \,:=\, \min(0,u_\eps-\lambda_{\eps,+})\,\leq\, 0,
\end{gather*}
we multiply \eqref{eq:lambda-diff} by $(u_\eps-\lambda_{\eps,+})_-\phi^2$,
and integrate over $\Omega_2$. We then deduce that
\begin{align}
  \int\limits_{\Omega_2} \eps|\nabla(u_\eps-\lambda_{\eps,+})_-|^2\phi^2 \,\leq\,&
  -\int_{\Omega_2} \eps (u_\eps-\lambda_{\eps,+})_-
  \nabla u_\eps\cdot2\phi\nabla\phi \notag\\ 
  &-\int_{\Omega_2}  \frac{1}{\eps}\Big(W'(u_\eps)-
  W'(\lambda_{\eps,+})\Big)(u_\eps-\lambda_{\eps,+})_-\phi^2\notag \\    
  \leq\, &\int_{\Omega_2}
  \eps |\nabla(u_\eps-\lambda_{\eps,+})_-|^2\phi^2 +
  \int_{\Omega_2}\eps |\nabla\phi|^2(u_\eps-\lambda_{\eps,+})_-^2\notag\\
  &-\frac{1}{2\eps}W''(1)\int_{\Omega_2}(u_\eps-\lambda_{\eps,+})_-^2\phi^2,
  \label{eq:est-u-lambda1}   
\end{align}
where in the last line we have used \eqref{eq:def-c0}, \eqref{eq:u-out}.
We therefore obtain that
\begin{gather}
   \frac{1}{2}W''(1)\int_{\Omega_1} (u_\eps-\lambda_{\eps,+})_-^2 \,\leq\,
   \eps^2\|\nabla\phi\|_{L^\infty(\Omega_2)}^2
   \int_{\Omega_2}(u_\eps-\lambda_{\eps,+})_-^2.
   \label{eq:est-u-lambda2}
\end{gather}
Choosing now bounded domains $\Omega_j$, $j=1,...,k+1$, such that
\begin{gather*}
  U\,=\,\Omega_1\,\subset\subset\,\Omega_2\,
  \subset\subset\,...\,\subset\subset\,\Omega_{k+1}\,=\, \Omega''
\end{gather*}
and iterating the procedure above we deduce that
\begin{gather}
   \int_{U} (u_\eps-\lambda_{\eps,+})_-^2 \,\leq\,
   C(k,U,\Omega'')\eps^{2k}\int_{\Omega''}(u_\eps-\lambda_{\eps,+})_-^2\,
   \leq\,C(k,U,\Omega'')\eps^{2k}.  
   \label{eq:est-u-lambda3}
\end{gather}
Assume now that for a $x_1\in U$
\begin{gather}
  u_\eps(x_1) \,<\, \lambda_{\eps,+} - \smallc\eps^{k}, \label{ass:u-lambda}
\end{gather}
choose $r=r(U,\Omega'')$ such that $B_{r\eps}(x_1)\subset\Omega''$,
and consider the scaled functions $\tilde{u},\tilde{f}:B_r^n(0)\to\R$,
\begin{gather*}
  \tilde{u}(x)\,:=\, u_\eps(x_1+\eps x), \qquad
  \tilde{f}(x)\,:=\, f_\eps(x_1+\eps x).
\end{gather*}
Then $\tilde{u},\tilde{f}$ satisfy the equation
\begin{gather*}
  -\Delta \tilde{u} \,=\, -W'(\tilde{u}) + \eps\tilde{f}.
\end{gather*}
Since the right-hand side is uniformly bounded we deduce that
$\tilde{u}\in W^{2,q}(B_r(0))$ for all $1\leq q<\infty$ and by the Sobolev
inequality that 
\begin{gather*}
  |\nabla\tilde{u}|\,\leq\, c_2(n,r)\quad\text{ on }B_r(0),\\
  |\nabla u_\eps|\,\leq\, \frac{c_2(n,r)}{\eps}\quad\text{ on }B_{\eps r}(x_1).
\end{gather*}
Thus \eqref{ass:u-lambda} gives us for all $x\in B_{r\eps^{k+1}}(x_1)$
\begin{gather*}
  u_\eps(x)\,=\, u_\eps(x_1) + u_\eps(x)-u_\eps(x_1)\,<\, \lambda_{\eps,+}
  -(\smallc-c_2(n,r))\eps^{k}
\end{gather*}
and we compute that for $\smallc>c_2(n,r)$
\begin{gather*}
  \int_{U} (u_\eps-\lambda_{\eps,+})_-^2 \,\geq\,
  \int_{B_{r\eps^{k+1}}(x_1)} (\smallc-c_2(n,r))^2\eps^k\,\geq\,
  (\smallc-c_2(n,r))^2\omega_n\eps^{k+n(k+1)}r^n. 
\end{gather*}
On the other hand, by \eqref{eq:est-u-lambda3} with $2k$ replaced by
$k+n(k+1)$ we obtain that
\begin{gather*}
   \int_{U} (u_\eps-\lambda_{\eps,+})_-^2 \,
   \leq\,\tilde{C}(k,U,\Omega'')\eps^{k+n(k+1)},  
\end{gather*}
which gives a contradiction for all $\smallc=\smallc(n,k,U,\Omega'')$
sufficiently large. 

To prove \eqref{eq:asym-u-eps-2} we first observe that
\begin{gather*}
  -\eps\Delta (-1+c_0) + \frac{1}{\eps}W'(-1+c_0) \,\geq\, \frac{8}{9}f(y_1)
\end{gather*}
for all $\eps>0$ sufficiently small. Since the minimum of two
supersolutions is a supersolution we deduce that
\begin{gather*}
  \tilde{u}_\eps\,:=\, \min(u_\eps, -1+c_0)
\end{gather*}
satisfies $\tilde{u}_\eps\in W^{1,2}(\Omega')$ and
\begin{gather*}
  -\eps\Delta \tilde{u}_\eps + \frac{1}{\eps}W'(\tilde{u}_\eps) \,\geq\,
   \frac{8}{9}f(y_1),\\
   \tilde{u}_\eps\,\leq\, -1+c_0
\end{gather*}
Then we can prove by the same arguments as for (1) that
\begin{gather*}
  \tilde{u}_\eps\,\geq\,
  \lambda_{\eps,-} - C_k\eps^k \quad\text{ in }U.
\end{gather*}
Since $\lambda_{\eps,-} - C_k\eps^k < -1+c_0$ for $\eps<\eps_0(k)$ this
proves \eqref{eq:asym-u-eps-2}.
\end{proof}
We will employ a comparison principle on the cylinder
$B_{\varrho_2}(y_1)\times (-4\varrho_0,4\varrho_0)$.
We first control the difference
$u_\eps-v_\eps$ on the top and the bottom, starting with the following
lemma.
\begin{lemma}
\label{lem:asym-val}
Let $\lambda_{\eps,+}, \lambda_{\eps,-}$ be the positive and
negative solution of \eqref{eq:def-lambda-eps} 
and let $\beta_{\eps,\pm}$ denote the values of $v_\eps$ `away' from
$M$, 
\begin{align}
  \beta_{\eps,+}\,&:=\, \phi_0(\eps^{-1}\delta(\eps)) +
  \eps\phi_1(\eps^{-1}\delta(\eps))\frac{2}{3\sigma}f(y_1),
  \label{eq:def-eps-beta-+}\\
  \beta_{\eps,-}\,&:=\, \phi_0(-\eps^{-1}\delta(\eps)) +
  \eps\phi_1(-\eps^{-1}\delta(\eps))\frac{2}{3\sigma}f(y_1).
  \label{eq:def-eps-beta--}
\end{align}
Then there exists $\gamma>0$, $\eps_0>0$ such that for all $0<\eps<\eps_0$
\begin{align}
  \lambda_{\eps,+} - \beta_{\eps,+} \,&\geq\, \gamma\eps,
  \label{eq:asym-val+}\\ 
  \lambda_{\eps,-} - \beta_{\eps,-} \,&\geq\,
  \gamma\eps. \label{eq:asym-val-} 
\end{align}
\end{lemma}
\begin{proof}
By a Taylor approximation
\begin{align*}
  W'(\beta_{\eps,\pm})\,&=\, W'\big(\phi_0(-\eps^{-1}\delta(\eps))\big)
  +  W''\big(\phi_0(-\eps^{-1}\delta(\eps))\big)\cdot
  \eps\phi_1(-\eps^{-1}\delta(\eps))\frac{2}{3\sigma}f(y_1)\\
  &\qquad +O(\eps^2)\\
  &=\, \frac{2}{3}f(y_1)\eps + o(\eps),
\end{align*}
where we have used that
\begin{gather*}
  W'\big(\phi_0(-\eps^{-1}\delta(\eps))\big)\,=\, o(\eps)
\end{gather*}
by \eqref{eq:phi}, \eqref{eq:cond-delta-3} and that
\begin{gather*}
  W''\big(\phi_0(-\eps^{-1}\delta(\eps))\big)
  \phi_1(-\eps^{-1}\delta(\eps))\frac{1}{\sigma} 
  \,\to\,0\quad\text{ as }\eps\to 0
\end{gather*}
by \eqref{eq:infty-phi_1}, \eqref{eq:cond-delta}.

We therefore deduce that
\begin{align}
  (\lambda_{\eps,\pm}-\beta_{\eps,\pm})\int_0^1 W''\big(s\lambda_{\eps,\pm}
  +(1-s)\beta_{\eps,\pm}\big)\,ds\,&=\,
  W'(\lambda_{\eps,\pm})-W'(\beta_{\eps,\pm})\notag\\
  &=\,
  \frac{1}{6}f(y_1)\eps +o(\eps). \label{eq:asym-prelim}
\end{align}
Since $\lambda_{\eps,-},\beta_{\eps,-}$ converge to $-1$ as
$\eps\to 0$ and since $W''(-1)>0$ we deduce from \eqref{eq:asym-prelim}
that \eqref{eq:asym-val-} holds for $\gamma>0$ and $\eps>0$
sufficiently small. Analogously
we obtain \eqref{eq:asym-val+}.
\end{proof}
\begin{proposition}\label{prop:comp}
For all $\eps>0$ sufficiently small we obtain that
\begin{gather}
\label{eq:comp}
  u_\eps  \,\geq\, v_\eps \quad\text{ in }
  \overline{B_{\varrho_2}(y_1)}\times [-4\varrho_0,4\varrho_0].
\end{gather}
\end{proposition}
\begin{proof}
Let us define the sets
\begin{gather*}
  U\,:=\, B_{\varrho_2}(y_1)\times (-4\varrho_0,4\varrho_0),\qquad
  \Omega'\,:=\, B_{\varrho}(y_1)\times (-5\varrho_0,5\varrho_0).
\end{gather*}
Consider for $s>0$ the \emph{shifted functions} $v_\eps^{(s)}$,
\begin{gather*}
  v_\eps^{(s)}(y,t)\,:=\, v_\eps(y,t-s)\quad\text{ for }(y,t)\in\Omega',
\end{gather*}
and the function
\begin{gather*}
  \Phi(s)\,:=\, \min_{\overline{U}} (u_\eps -v_\eps^{(s)}).
\end{gather*}
Assume now $\eps<\eps_1$, where we choose $\eps_1>0$ below, and that
\eqref{eq:comp} is not satisfied, hence
\begin{gather}
  \Phi(0)\,<\, 0. \label{ass:Phi}
\end{gather}
The definition of $v_\eps$ in \eqref{eq:def-beta}, \eqref{eq:def-v-eps}
implies that we can choose $s_0>0$, $s_0=s_0(\varrho_0)$ such that for
all $\eps>0$ sufficiently small
\begin{gather}
  v_\eps^{(s_0)}\,=\, \beta_{\eps,-}\quad\text{ in
  }\Omega'. \label{eq:v-eps-s} 
\end{gather}
Applying then Lemma \ref{lem:bdry-1} with $U,\Omega'$ as above and $k=2$
we deduce that
\begin{gather}
  u_\eps\,\geq\,\lambda_{\eps,-} -O(\eps^2)\quad\text{ in
  }U. \label{eq:u-eps-s} 
\end{gather}
Therefore \eqref{eq:asym-val-} and \eqref{eq:v-eps-s},
\eqref{eq:u-eps-s} imply that
\begin{gather*}
  u_\eps - v_\eps^{(s_0)}\,\geq\, \lambda_{\eps,-}
  -O(\eps^2)-\beta_{\eps,-} \,\geq\, \gamma\eps -O(\eps^2)
\end{gather*}
and we deduce that
\begin{gather}
  \Phi(s_0)\,>\, \frac{\gamma}{2}\eps \label{eq:Phi-s-0}
\end{gather}
for all $\eps<\eps_1$ and $\eps_1>0$ chosen suitably small.

Since $\Phi$ is continuous \eqref{ass:Phi}, \eqref{eq:Phi-s-0} imply the
existence of $s_*>0$ and $x_*\in\overline{U}$ such that
\begin{gather}
   0\,=\,(u_\eps-v_\eps^{(s_*)})(x_*)\,=\, \min_{\overline{U}} (u_\eps
   -v_\eps^{(s_*)}). \label{eq:min-0}
\end{gather}
We first prove that $x_*\in U$. With  this aim we consider the different
parts of $\partial U$.
\begin{enumerate}[leftmargin=*]
\item
For $\eps_1>0$ chosen suitably small we have
\begin{gather*}
  v_\eps^{(s_*)}\,=\, \beta_{\eps,-}\quad\text{ in
  }B_{\varrho_2}(y_1)\times (-5\varrho_0,-4\varrho_0)
\end{gather*}
and therefore, 
by Lemma \ref{lem:bdry-1} and Lemma \ref{lem:asym-val},
\begin{gather*}
  u_\eps\,\geq\, v_\eps^{(s_*)} + \frac{\gamma}{2}\eps \quad\text{ in
  }B_{\varrho_2}(y_1)\times (-5\varrho_0,-4\varrho_0)
\end{gather*}
for $\eps<\eps_1$ and $\eps_1>0$ suitably small, see the argument
above. 
This shows that
\begin{gather}
  x_*\,\not\in \,
  \overline{B_{\varrho_2}(y_1)}\times\{-4\varrho_0\}. \label{eq:bdry1} 
\end{gather}
\item
For $\eps_1>0$ sufficiently small we obtain from Lemma \ref{lem:To} and
Lemma \ref{lem:bdry-1}, applied to
\begin{gather*}
  U_+\,=\, B_{\varrho_2}(y_1)\times (2\varrho_0,4\varrho_0),\quad
  \Omega''_+\,=\, B_{\varrho}(y_1)\times (\varrho_0,5\varrho_0), k=2 
\end{gather*}
that
\begin{gather*}
  u_\eps\,\geq\, \lambda_{\eps,+}-O(\eps^2)\quad\text{ in }U_+.
\end{gather*}
By \eqref{eq:asym-val+} this implies that
\begin{gather*}
  u_\eps\,\geq\,\beta_{\eps,+} +\gamma_\eps -O(\eps^2)\,\geq\,
  \beta_{\eps,+} +\frac{\gamma}{2}\eps\,\geq\, v_\eps +\frac{\gamma}{2}
\end{gather*}
holds in $U_+$, hence
\begin{gather}
  x_*\,\not\in\,
  \overline{B_{\varrho_2}(y_1)}\times\{4\varrho_0\}. \label{eq:bdry2} 
\end{gather}
\item
By \eqref{eq:psi-phi-bdry} there exists bounded domains
$U_1,\Omega_1,U_2$ such that
\begin{align}
  & U_1\,\subset\subset\,\Omega_1\,\subset\subset
  B_{\varrho}(y_1)\setminus\overline{B_{\varrho_3}(y_1)}\times
  (-5\varrho_0,5\varrho_0) \\
  & \Omega_1\,\subset\subset\, \{t>\varphi_+(y)\},\\
  & U_2\,\subset\subset\, \{t<\psi(y)\},\\
  & \partial B_{\varrho_2}(y_1)\times
  [-4\varrho_0,4\varrho_0]\,\subset\, U_1\cup U_2. \label{eq:decomp}
\end{align}
By similar arguments as above we first prove that $u_\eps>v_\eps^{(s_*)}$ in
$U_1$. First we obtain from Lemma \ref{lem:To} that
\begin{gather*}
  u_\eps\,\geq\, 1-c_0\quad\text{ in }\Omega_1
\end{gather*}
and applying Lemma \ref{lem:bdry-1} and Lemma \ref{lem:asym-val} with
$k=2$ we deduce that in $U_1$
\begin{gather}
  u_\eps \,\geq\, \lambda_{\eps,+} -O(\eps^2)\,>\,
  v_\eps^{(s_*)}. \label{eq:U1-part} 
\end{gather}
Since $s_*>0$ we obtain that for $\eps<\eps_1$,
where $\eps_1=\eps_1(U_2,\varrho_2)$ is chosen sufficiently small,
\begin{gather}
  v_\eps^{(s_*)}\,\leq\, \beta_{\eps,-} +o(\eps)\,\leq\,
  \lambda_{\eps,-}-\gamma\eps +o(\eps) \,<\, u_\eps \label{eq:U2-part}
\end{gather}
in $U_2$.
By \eqref{eq:decomp}, \eqref{eq:U1-part}, \eqref{eq:U2-part} we deduce
that 
\begin{gather}
  u_\eps\,>\,v_\eps^{(s_*)}\quad\text{ on }\partial B_{\varrho_2}(y_1)\times
  (-4\varrho_0,4\varrho_0). \label{eq:bdry3}
\end{gather}
\end{enumerate}
By \eqref{eq:min-0}-\eqref{eq:bdry2}, \eqref{eq:bdry3} we get that
$u_\eps-v_\eps^{(s_*)}$ 
has an interior minimum with value zero at $x_*$.
Using \eqref{eq:GT-eps}, \eqref{eq:rhs-u}, and \eqref{eq:gt-v-eps-global}
we therefore deduce that
\begin{align*}
  0\,&\geq\, -\eps\Delta
    (u_\eps-v_\eps^{(s_*)})(x_*)\\
  &\geq\, f_\eps(x_*) -\frac{7}{9}f(y_1) 
  - \frac{1}{\eps}W'(u_\eps(x_*)) +
  \frac{1}{\eps}W'(v_\eps^{(s_*)}(x_*)) \\
  &\geq\,\frac{1}{9}f(y_1).
\end{align*}
This finally gives a contradiction
and proves Proposition \ref{prop:comp-psi}.
\end{proof}
%=======================================================
% applications
%=======================================================
\section{Applications}
\label{sec:applications}

\subsection{Stationary solutions for the Cahn-Hilliard equation}
\label{subsec:CH}

Consider the Cahn-Hilliard equation \cite{CH}
\begin{gather}
  \frac{\partial u_{\varepsilon}}{\partial t}=\Delta
  f_{\varepsilon},\qquad f_\eps\,=\, -\eps\Delta u_\eps +
  \frac{1}{\eps}W'(u_\eps) \label{eq:CH}
\end{gather}
with Neumann boundary conditions for $u_{\varepsilon}$ and
$f_{\varepsilon}$.

The stationary solutions of \eqref{eq:CH} are those with constant
$f_{\varepsilon}$,
\begin{gather}
  -\eps\Delta u_\eps + \frac{1}{\eps}W'(u_\eps)\,=\,
  \lambda_\eps,\quad \lambda_\eps\in\R. \label{eq:stat-CH}
\end{gather}
This is also the Euler--Lagrange equation of the volume constrained
minimization problem for the Cahn--Hilliard energy \eqref{eq:def-E},
\begin{gather}
  \inf \big\{E_\eps(u)\,:\, u\in H^{1,2}(\Omega), \int_\Omega u\,=\,
  m\}.   \label{eq:E-min-vol}
\end{gather}
To better understand stationary solutions of the
Cahn--Hilliard equation the asymptotic of \eqref{eq:stat-CH} is
analyzed in a couple of papers. The behavior of (locally) energy
minimizing solutions of \eqref{eq:stat-CH} is well understood
\cite{Mod,Ste,LM,HT}. In this case sequence $(u_\eps)_{\eps>0}$ with uniformly bounded
energy converge to a constant-mean curvature hypersurface with
single-multiplicity. This hypersurface is smooth except for a
closed set of dimension at most $n-8$.

Solutions of the  Cahn--Hilliard equation are observed to undergo
pattern similar to unstable equilibria \cite{GM88} and the behavior of
general stationary points is another question of interest. However, this
situation is more difficult due to the possibility of higher-multiplicity
surfaces in the limit. Hutchinson--Tonegawa \cite{HT} showed that the limit is
given by an integer-rectifiable varifold and that the weak mean
curvature exists. However this mean curvature is only \emph{locally}
constant, the constant depending 
on the multiplicity. A higher multiplicity is also an
obstacle to use Allard's regularity theory \cite{All} and to obtain the
smoothness of the limit.

As a corollary of Theorem \ref{the:claim} we can improve the previous
results. 
\begin{theorem}
\label{the:stat-CH}
Consider a sequence $(u_{\varepsilon})_{\eps>0}\subset H^{1,2}(\Omega)$
with a fixed volume constraint $\int_{\Omega}u_{\varepsilon}=m$, and a
sequence $(\lambda_\eps)_{\eps>0}\subset\R$ such that \eqref{eq:stat-CH}
is satisfied.
Suppose further that
\begin{gather}
  E_{\varepsilon}(u_{\varepsilon}) \,\leq\, \Lambda\quad\text{ for all
  }\eps>0. \label{eq:ass-CH-E}
\end{gather}

Then there exists a subsequence $\eps\to 0$, a function $u\in
BV(\Omega,\{-1,1\})$, and $\lambda:=\lim_{\eps\to 0}\lambda_\eps$, such that
$u_\eps\to u$ in $L^1(\Omega)$. The phase boundary $\partial^*\{u=1\}\cap\Omega$
has constant mean curvature $\sigma H=\lambda$. 

In the case that
$\lambda\neq 0$  the phase boundary is up to a $\Ha^{n-1}$-nullset a
smooth hypersurface.
The energy measures $\mu_\eps$ as defined in
\eqref{eq:def-mu-eps} converge to a measure $\mu$ that is up to the
factor $2\sigma$ integer-rectifiable, has constant mean curvature
$\sigma H=\lambda$ and multiplicity one $\Ha^{n-1}$-almost everywhere on
$\partial^*\{u=1\}$. Moreover `hidden boundaries' can only occur in one
phase and have zero mean curvature:
\begin{gather}
  \mu(K)\,=\, 0 \quad
  \begin{cases}
    \text{ for all compact sets }K\subset \{u=1\}^o &\text{ if
    }\lambda > 0,\\
    \text{ for all compact sets }K\subset \{u=-1\}^o &\text{ if
    }\lambda < 0,
  \end{cases} \label{eq:hidden}\\
  H\,=\, 0 \quad\Ha^{n-1}\text{-almost everywhere on }\spt(\mu)\setminus
  \partial^*\{u=1\}. \label{eq:hidden2}
\end{gather}
\end{theorem}
\begin{proof}
It follows from \cite[Lemma 3.4]{Che} that $|\lambda_{\varepsilon} |\leq
c(m,\Lambda)$ and we may choose a subsequence such that
$\lambda=\lim_{\eps\to 0} \lambda_{\varepsilon}$ exists.

We therefore can apply Theorem \ref{the:claim} and obtain that there
exists a subsequence $\eps\to 0$ and limits $u,\mu$ of $u_\eps,\mu_\eps$.
Moreover, $u\in BV(\Omega,\{-1,1\})$ and
$(2\sigma)^{-1}\mu$ is an integer-rectifiable varifold with weak mean
curvature
\begin{gather*}
  \sigma H\,=\, \lambda\quad\Ha^{n-1}-\text{almost everywhere on
  }\partial^*\{u=1\}.
\end{gather*}
Next it follows from Theorem \ref{the:claim} that $\lambda=0$ on the
parts of $\partial^*\{u=1\}\cap\Omega$ with odd multiplicity larger than $1$, which
shows that in the case $\lambda\neq 0$ the phase boundary is given as a
constant curvature varifold with unit multiplicity. By Allard's
regularity theory \cite{All} we conclude the smoothness of the phase
boundary.
Finally \eqref{eq:hidden}, \eqref{eq:hidden2} follow from
\eqref{eq:claim-even}. 
\end{proof}
In general dimension we can not insure good regularity of the hidden
boundaries, due to the lack of regularity theory for general stationary
integral varifold. Only for $n=2$ we can conclude that ${\rm spt}\,
\mu\cap K$ is given by straight line segments with possible 
junction points for all compact sets $K\subset\Omega\setminus
\partial^*\{u=1\}$ \cite{AA}.

\subsection{Critical points of the Ohta-Kawasaki functional}
\label{subsec:OK}
The micro-phase separation of block copolymers exhibits the formation of
complex patterns. Ohta--Kawasaki \cite{OK86} and later Bahiana--Oono
\cite{BO90} used a 
phase-field like approach and proposed a free energy that is
after a suitable rescaling given by
\begin{gather}
  F_{\varepsilon}(u) \,=\, \int_\Omega \Big(\frac{\varepsilon}{2}
  |\nabla u|^2+\frac{1}{\varepsilon}W(u)+
  \frac{1}{2}|\nabla v|^2\Big) dx, \label{eq:F-OK}
\end{gather}
where $v(\cdot)=v[u](\cdot)$ is the solution of
\begin{gather}
  -\Delta v \,=\, u- \frac{1}{|\Omega|}\int_\Omega u\quad\text{ in
  }\Omega,
  \qquad \nabla v\cdot\nu_\Omega\,=\,0 \quad\text{ on
  }\partial\Omega. \label{eq:OK-v}
\end{gather}
The functional $F_\eps$ extends the Cahn--Hilliard energy by a non-local
term that describes long-range interactions between chains of
macromolecules. For a derivation of this model by a density-functional
approach see \cite{CR03}. The set of (local) minimizers of
$F_\eps$ is extremely rich and \eqref{eq:F-OK} has drawn quite some
attention \cite{Cho01,ONIM,RW03,RW06}. The Gamma-limit of
$F_\eps$ as $\eps\to 0$ and the convergence of the
corresponding $H^{-1}$ gradient-flow that was proposed by Nishiura and
Ohnishi \cite{NO95} are also well-studied \cite{RW00,FH01,Hen01}. 

Critical points of $F_{\varepsilon}$ under a volume-constraint satisfy
the Euler--Lagrange equation 
\begin{gather}
  -\varepsilon\Delta u+\frac{1}{\varepsilon}W'(u) +v[u]\,=\,
  \lambda, \label{eq:EL-OK}
\end{gather}
where $\lambda\in\R$ is a Lagrange-multiplier.

As a corollary of our results we obtain the convergence of stationary
points of $F_\eps$.
\begin{theorem}
\label{the:OK}
Assume that we have a sequence
$(u_{\varepsilon})_{\eps>0}$ such that
\begin{gather}
  F_{\varepsilon}(u_{\varepsilon})+\|u_{\varepsilon}\|_{L^{\infty}(\Omega)}
  \,\leq\,
  \Lambda\quad\text{ for all }\eps>0 \label{eq:bound-OK}
\end{gather}
and such that \eqref{eq:EL-OK} holds for Lagrange multipliers
$\lambda_\eps\in\R$ and the solutions $v_\eps$ of
\begin{gather}
  -\Delta v_\eps \,=\, u_\eps- \frac{1}{|\Omega|}\int_\Omega
  u_\eps\quad\text{ in 
  }\Omega,
  \qquad \nabla v_\eps\cdot\nu_\Omega\,=\,0 \quad\text{ on
  }\partial\Omega. \label{eq:OK-v-eps}
\end{gather}

Then there exists a subsequence $\eps\to 0$, a number $\lambda\in\R$,
and a function $u\in BV(\Omega,\{-1,1\})$ such that
$\lambda=\lim_{\eps\to 0}\lambda_\eps$ and
$u_\eps\to u$ in $L^1(\Omega)$. Moreover $v_\eps\to v$ in
$C^{1,\alpha}(\bar{\Omega})$ for all $0<\alpha<1$ and $v$ solves
\eqref{eq:OK-v}.

The energy measures $\mu_\eps$ as defined in
\eqref{eq:def-mu-eps} converge to a measure $\mu$ that is up to the
factor $2\sigma$ integer-rectifiable and has a weak mean curvature that
satisfies
\begin{gather}
  \sigma H\,=\,
  \begin{cases}
    -v + \lambda\quad&\Ha^{n-1}-\text{almost everywhere on }
    \partial^*\{u=1\}, \\
    0 &\Ha^{n-1}-\text{almost everywhere on }\spt(\mu)\setminus
    \partial^*\{u=1\}.
  \end{cases} \label{eq:GT-OK}
\end{gather}
Finally $\partial^*\{u=1\}$ has multiplicity one $\Ha^{n-1}$-almost everywhere
in the set $\{v\neq\lambda\}$ and this part of the phase boundary is a
$C^{3,\alpha}$-surface for all $\alpha<1$, except for a set of 
${\mathcal H}^{n-1}$-measure zero.
\end{theorem}
\begin{proof}
By standard elliptic theory we obtain from \eqref{eq:bound-OK},
\eqref{eq:OK-v-eps} that $v_\eps$ 
is uniformly bounded in $W^{2,p}(\Omega)$ for all $1\leq
p<\infty$. Therefore Theorem \ref{the:claim} applies and we can repeat
the arguments of the proof of Theorem \ref{the:stat-CH}. We omit the
details here.
\end{proof}
If $u_{\varepsilon}$ has in addition a local energy
minimizing property for $F_{\varepsilon}$ we can draw stronger
conclusions: Then $\mu$ has multiplicity one $\mu$-almost everywhere and
is $C^{3,\alpha}$-smooth, see the arguments in \cite{HT}.
%
%==================================
% appendix
%==================================
\begin{appendix}
\section{A generalization of mean curvature to general phase boundaries}
\label{app:defs}
\begin{proposition}[see {\cite[Proposition 3.1]{Ro04}}]\label{prop-def-h}
Let\/ $\Omega\subset\R^n$ be open, $E\subset\Omega$, and $\Chi_E\in \BV(\Omega)$. Assume that there are two
integral\/ $(n-1)$-varifolds $\mu_1,\mu_2$ on\/ $\Omega$ such that for $i=1,2$ the following hold:
\begin{gather}
\partial^*E \,\subset\, \spt (\mu_i),\label{ass-def-h1}\\
\mu_i\text{ has locally bounded first variation with mean curvature vector }\vec{H}_{\mu_i},\label{ass-def-h2}\\
\vec{H}_{\mu_i}\,\in\, \Lp^s_{\loc}(\mu_i),\,s>n-1,\,s\geq 2\label{ass-def-h3}.
\end{gather}
Then
\[
\vec{H}_{\mu_1}|_{\partial^*E} = \vec{H}_{\mu_2}|_{\partial^*E}
\]
is satisfied\/ $\mathcal{H}^{n-1}$-almost everywhere on $\partial^*E$.
\end{proposition}

This proposition justifies the following definition.

\begin{definition}\label{def-h}
Let $E\subset\Omega$ and $\Chi_E\in \BV(\Omega)$, and assume that there exists an
integral\/ $(n-1)$-varifold $\mu$ on\/ $\Omega$ satisfying\/ {\rm(\ref{ass-def-h1})--(\ref{ass-def-h3})}.
Then we call
\[
\vec{H} := \vec{H}_{\mu}|_{\partial^*E}
\]
the generalized mean curvature vector of $\partial^*E$ and define a scalar mean curvature by
\[
H := \vec{H}\cdot \frac{\nabla\Chi}{|\nabla\Chi|}\text{ on }\partial^*E.
\]
\end{definition}
\end{appendix}
%%%%%%%%%%%%%%%%%%%%%%%%
% bibliography
%%%%%%%%%%%%%%%%%%%%%%%%

%\bibliography{gt-lib}
%\bibliographystyle{plain}
%

\end{document}